\documentclass[10pt]{amsart}
\usepackage{amsfonts}
\usepackage{indentfirst,bm,amsmath}
\usepackage{amsthm,amssymb}
\usepackage{mathrsfs}
\usepackage{float}
\usepackage{graphicx}
\usepackage{subfig}
\usepackage{amsmath}
\usepackage{amssymb}
\usepackage{hyperref}
\usepackage[numbers,sort&compress]{natbib}

\newtheorem{theorem}{Theorem}[section]

\newtheorem{lemma}[theorem]{Lemma}
\newtheorem{proposition}[theorem]{Proposition}
\newtheorem{definition}[theorem]{Definition}
\newtheorem{remark}[theorem]{Remark}
\newtheorem{example}[theorem]{Example}

\begin{document}
	\title[]{Nevanlinna theory for generalized tropical meromorphic functions and applications to ultra-discrete equations}
	
	\author[T. B. Cao, R. Korhonen \and  Y. T. Wang]{ Tingbin Cao, Risto Korhonen \and Yuting Wang}
	
    \address[Tingbin Cao]{Department of Mathematics, Nanchang University, Nanchang city, Jiangxi 330031, P. R. China}
    \address[Risto Korhonen]{Department of Physics and Mathematics, University of Eastern Finland, P.O. Box 111 FI-80101 Joensuu, Finland}
    \address[Yuting Wang]{Department of Mathematics, Nanchang University, Nanchang city, Jiangxi 330031, P. R. China}
	\email{tbcao@ncu.edu.cn, risto.korhonen@uef.fi, ytwang@email.ncu.edu.cn}
	
	\thanks{The first author is supported by the National Natural Science Foundation of China (\#12571082), the Jiangxi Natural Science Foundation (\#20232ACB201005) and the Shandong Natural Science Foundation (\#ZR2024MA024). The third author is supported by China Scholarship Council (\#202606820022).}
	\date{}
	
	
	
	\begin{abstract} The main purpose of this paper is to establish the tropical Nevanlinna theory to piecewise linear real functions with arbitrary first-kind discontinuities (called generalized tropical meromorphic functions). The key point is that we introduce the directional jump counting functions $J^+(r,f)$ and $J^-(r,f)$ for the case of arbitrary first-kind discontinuities to refine the single function $J(r,f)$ in the previous theory. We first obtain the new tropical version of the Poisson-Jensen formula, the first main theorem. Under the assumption that \(f\) is slope-admissible, we show that \(T(r,f)\) is nondecreasing and establish a logarithmic derivative lemma.  These results unify the classical continuous case due to Halburd and Southall and the one-sided continuous case due to Halonen, Korhonen and Filipuk.\par

As applications, we study several first-order and second-order ultra-discrete equations, showing that certain equations admit infinitely many generalized tropical meromorphic solutions, and that under suitable conditions infinite-order solutions also exist. Finally, using the logarithmic derivative lemma obtained above, we establish a tropical Clunie lemma.

	\end{abstract}
	\maketitle
	\tableofcontents

\section{Introduction and main results} In recent over ten years, Nevanlinna theory in tropical geometry has been a very active and popular research topic. In 2009, Halburd and Southall \cite{Halburd09} first established tropical version of the Poisson-Jensen formula and logarithmic derivative lemma  for continuous piecewise linear real functions with integer slopes (called tropical meromorphic functions) in the tropical semifield  $\mathbb{T} = \mathbb{R} \cup \{-\infty\},$ where the addition and multiplication are defined, respectively, by
\[
x\oplus y:=\max(x,y)
\]
and
\[
x\otimes y:=x + y.
\] In addition, the tropical division and exponentiation are defined as \(x\oslash y:=\frac{x}{y}\oslash=x - y\) and \(x^{\otimes \alpha}:=\alpha x\) for \(\alpha\in\mathbb{R}\), respectively. The additive identity element is \(0_{\mathbb{T}}=-\infty\), and the multiplicative identity element is $1_{\mathbb{T}}=0$. For a more comprehensive account of tropical geometry, see references \cite{mik, Maclagan-15, Mik2006, Mik2020, Speyer09}.

Subsequently,  tropical Nevanlinna theory has developed rapidly, and many important works have emerged. In 2011, Laine and Tohge \cite{Laine11} extended the theory to tropical meromorphic functions whose slopes may be arbitrary real numbers, under the growth condition that the hyper-order is strictly less than one. In 2016, Korhonen and Tohge \cite{Risto16} considered the degeneracy of the set of tropical linear combinations of $\{f_0,\dots,f_n\}$, and developed tropical Nevanlinna theory for tropical holomorphic curves into $n$-dimensional tropical projective space. In 2025, Cao and Zheng \cite{Cao-Zheng25} further studied tropical hypersurfaces under a minimal hypertype growth condition. Very recently, Cao and Peng \cite{Cao-peng} successfully extended these results to several tropical variables. Korhonen and Tan \cite{Risto26} considered tropical Nevanlinna theory from piecewise linear functions to piecewise polynomial functions. Furthermore, Wang and Cao \cite{W-C26} considers the case where the tropical hyperplane is in general position. This new theory not only establishes an intersection between complex analysis and tropical geometry, but also provides powerful analytical tools for the study of ultra-discrete systems, tropical dynamical systems and even neural network theories (see \cite{Laine16, Orm2010, Zhang18}).\par

However, all these classical results are restricted to continuous tropical meromorphic functions. To break the continuity constraint, Halonen, Korhonen and Filipuk \cite{Halonen23} recently introduced a special class of piecewise linear functions which admitting jump discontinuities, and established a tropical Poisson-Jensen formula to such functions. Their framework imposes strict one-sided continuity requirements at every discontinuity point. For such functions \(f\), Halonen, Korhonen and Filipuk introduced the notions of positive and negative jumps as follows. If $x\Omega_f(x) > 0$ and $x \neq 0$, the point $x$ is said to be a \textit{positive jump}; if $x\Omega_f(x) < 0$ and $x \neq 0$, the point $x$ is said to be a \textit{negative jump}, where
\[
\Omega_f(x) := f(x^+) - f(x^-)= \lim_{\varepsilon\to 0_{+}} f(x+\varepsilon)-\lim_{\varepsilon\to 0_{+}} f(x-\varepsilon).
\]
At the origin, a discontinuity is called \textit{positive jump} if $f(0) = \min\{f(0^+), f(0^-)\}$, and \textit{negative jump} if $f(0) = \max\{f(0^+), f(0^-)\}$. The corresponding \textit{jump counting function} is defined as:
\[
J(r, f) = \frac{1}{2} \sum_{|\beta_\nu| \leq r} h_f(\beta_\nu),
\]
where $\beta_\nu$ are the negative jump points of $f$ and $h_f(x) := |\Omega_f(x)|= |f(x^+) - f(x^-)|.$\par

\begin{figure}[H]
  \centering
  \includegraphics[width=0.35\linewidth]{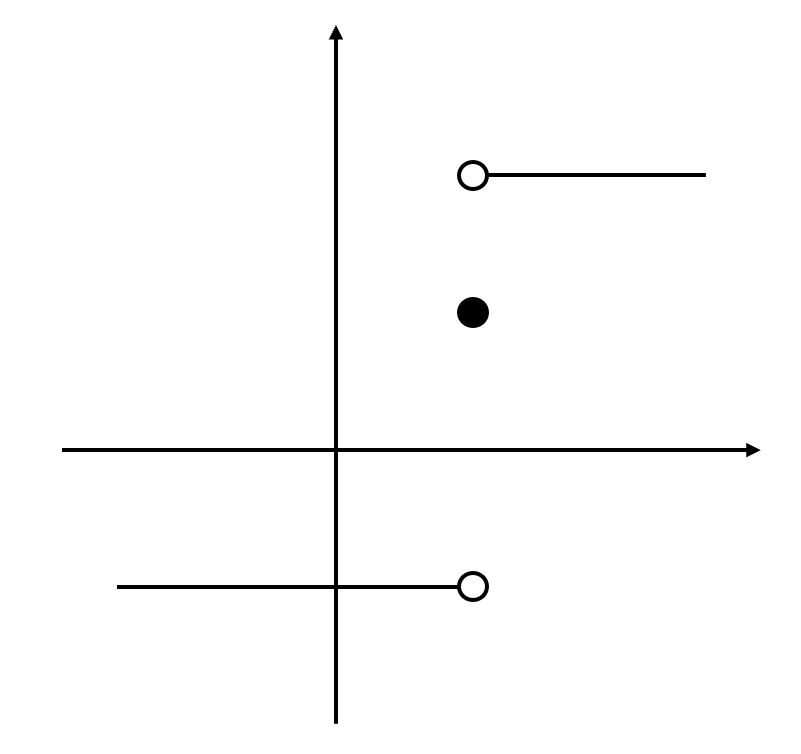}\\
  \caption{Function $f$}\label{1.jpg}
\end{figure}

However, we observe that this special class of one-sided continuous piecewise linear functions is not closed under the tropical operations. For example,  the point $x=1$ is a left-continuous and right-continuous of the functions
\[
f_1(x)=
\begin{cases}
-2, & x \leq 1, \\
2, & x> 1,
\end{cases}\quad \mbox{and}\quad
f_2(x)=
\begin{cases}
-1, & x < 1, \\
1, & x\geq 1,\end{cases}
\] respectively. However, $x=1$ is neither left-continuous nor right-continuous of $f=f_1\oplus f_2$ (see Figure \ref{1.jpg}). This means that the $f$ is not piecewise linear in the means of Halonen, Korhonen and Filipuk.

Motivated by this observation, we consider piecewise linear real functions with first-kind discontinuities as defined below (we call them generalized tropical meromorphic functions throughout this paper). This class of functions is closed under the tropical operations, a property that will be proved in Proposition \ref{L11}.

\begin{definition}[Generalized tropical meromorphic function]\label{D2.1} Let $f: \mathbb{R}\to\mathbb{R}.$ If there exist finitely or countably many partition points $\{x_k\}$ satisfying $$\cdots<x_{-2}<x_{-1}<x_0<x_1<x_2<\cdots,$$ such that for each $k,$ $f(x)=\alpha_k x+\beta_k$ for all $x\in(x_{k-1}, x_k)$, where $\alpha_k, \beta_k\in\mathbb{R},$ and $f$ is defined at every partition point $x_k$, with the partition points forming a locally finite set. Then $f$ is called a piecewise-linear function.
\end{definition}

For generalized tropical meromorphic functions $f$, we provide a more refined characterization of their right-direction jumps and left-direction jumps. Specifically, we say that $x$ is a \textit{right-positive jump} if $sgn(x^+)\Omega_f(x^+) > 0$ and a \textit{right-negative jump} if $sgn(x^+)\Omega_f(x^+) < 0$, where
\[
\Omega_f(x^+) := f(x^+) - f(x)= \lim_{\varepsilon\to 0_{+}} f(x+\varepsilon)- f(x).
\] Similarly, we say that $x$ is a \textit{left-positive jump} if $sgn(x^-)\Omega_f(x^-) > 0$, and a \textit{left-negative jump} if $sgn(x^-)\Omega_f(x^-) < 0$, where
\[
\Omega_f(x^-) := f(x) - f(x^-)=  f(x) -\lim_{\varepsilon\to 0_{+}} f(x-\varepsilon).
\] Clearly, we know that $\Omega_f(x^-)=\Omega_f(x^+)=0$ if and only if $f$ is continuous at $x$. It can be seen that, unlike the discussion in the sense of \cite{Halonen23}, our definition does not require a separate discussion of the case $x=0$.

The corresponding \textit{right height} and \textit{left height} at the jump point $x$ of $f$ are defined as
$$h_f(x^+) := |sgn(x^+)\Omega_f(x^+)|= |f(x^+) - f(x)|,$$
$$h_f(x^-) := |sgn(x^-)\Omega_f(x^-)|= |f(x) - f(x^-)|.$$ We then define the \textit{right jump counting function} of $f$ by
\[
J^+(r, f) = \frac{1}{2} \sum_{|\beta_\nu^+| \leq r} h_f(\beta_\nu^+),
\]
where $\beta_\nu$ are the right-negative jump points of $f$. The \textit{left jump counting function} of $f$ is defined by
\[
J^-(r, f) = \frac{1}{2} \sum_{|\zeta_t^-| \leq r} h_f(\zeta_t^-),
\]
where $\zeta_t$ are the left-negative jump points of $f$. For convenience, we denote the \textit{total jump counting function} by
\[
J(r,f):=J^+(r, f)+J^-(r, f).
\]

The characteristic function is defined by
\[
T(r,f):=m(r,f)+N(r,f)+J(r,f),
\]
where \(m(r,f)\) is the tropical proximity function and \(N(r,f)\) is the tropical counting function of poles. Our main results are as follows.

First, we prove a tropical Poisson--Jensen formula (Theorem \ref{T3.1}) which implies the following tropical Jensen formula.

\begin{theorem}[Jensen formula]\label{C3.2} Let $f$ be a generalized tropical meromorphic function in $\mathbb{R}$. Then
\begin{align*}
f(0)=\frac{1}{2}[f(r)+f(-r)]-N(r, 1_{\mathbb{T}}\oslash f) + N(r, f)-J(r, 1_{\mathbb{T}}\oslash f) + J(r, f)
\end{align*} or equivalently,
\[
T(r, f)=T(r, 1_{\mathbb{T}}\oslash f)+f(0)
\] for any $r>0$.
\end{theorem}

From this we derive the first main theorem:
\begin{theorem}\label{T1.4} Let $f$ and $a$ be generalized tropical meromorphic functions. Then
\begin{align}\label{eq7}
T\left(r,\frac{1_{\mathbb{T}}}{f\oplus a}\oslash\right)
&= T(r,f)-N(r,f)+N(r,f\oplus a)\notag\\
&\quad-J(r,f)+J(r,f\oplus a)-f(0)\oplus a(0)+\varepsilon_f(r,a)
\end{align}
for some quantity $\varepsilon_f(r,a)$ satisfying $0\le \varepsilon_f(r,a)\le m(r,a)$ for all $r\ge 0$.
\end{theorem}

Next, in Section 4, we give an example showing that $T(r,f)$ is not necessarily nondecreasing (see Example \ref{E4.1}); moreover, from this example, we can see that if right-derivative $f'(x^+)$ or left-derivative $f'(x^-)$ of generalized tropical meromorphic function $f$ is always positive, or always negative, then $T(r,f)$ may not be a non-decreasing function. Hence, we exclude these special cases and introduce the following definition.

\begin{definition}\label{D4.2} Let $f:\mathbb{R}\to\mathbb{R}$ is a  generalized tropical meromorphic function. We say that $f$ is slope-admissible if its right-derivative $f'(x^+)$ or left-derivative $f'(x^-)$ is neither always positive nor always negative for $x\in\mathbb{R}$.
\end{definition}

Interestingly, we find that $T(r, f)$ is nondecreasing when $f$ is slope-admissible (see Theorem \ref{T4.3}). Furthermore, we give an estimate of the proximity function of $f(x+c)\oslash f(x)$ ( Theorem \ref{T4.33} ) and establish the following logarithmic derivative lemma for slope-admissible generalized tropical meromorphic function with hyper-order strictly less than one, where the hyper-order of $f$ is defined by
\[
\rho_2(f):=\limsup_{r\to \infty}\frac{\log\log T(r,f)}{\log r}.
\] We note that due to Example \ref{E4.4},  $T(r,f)$ may be not necessarily convex, and thus the logarithmic derivative lemma can not be improved to the grow of minimal hypertype for a generalized tropical meromorphic function $f$ (i.e., \(\limsup_{r\to\infty}\frac{\log T_f(r)}{r}=0\))  using the method of Cao and Zheng \cite{Cao-Zheng25}.

\begin{lemma}[Logarithmic derivative lemma] Let $c\in\mathbb{R}\setminus\{0\}$. If $f$ is slope-admissible generalized tropical meromorphic function on $\mathbb{R}$ with
with hyper-order strictly less than one, then
\[m(r,f(x + c)\oslash f(x)) = o(T(r,f))\]
holds outside an exceptional set of finite logarithmic measure as $r$ tends to infinity.
\end{lemma}

In addition, we apply these results to first-order and second-order ultra-discrete equations. In particular, we consider the following ultra-discrete equation.
\[
y(x+1)\otimes y(x-1)=R\big(x,y(x)\big),
\] where $R$ is tropical rational in $x$ and $y$. We noted that all such equations admit infinitely many generalized tropical meromorphic solutions.

\begin{theorem}
Let $P(y)=\max\{a_0,a_1+y,\dots,a_p+py\}$ and $Q(y)=\max\{b_0,b_1+y,\dots,b_q+qy\}$ be two tropical polynomials with no common roots and neither $a_p$ nor $b_q$ is $-\infty$. If $|p-q|>2$, then the equation
\begin{equation*}
y(x+1)+y(x-1)=P(y(x))-Q(y(x))
\end{equation*}
admits infinitely many infinite-order generalized tropical meromorphic solutions.
\end{theorem}

Finally, by means of the logarithmic derivative lemma, we obtain a tropical Clunie-type lemma (Theorem \ref{T5.4}) for slope-admissible piecewise-linear functions.

This paper is organized as follows. In the section 2, we first briefly introduce  some basic definitions and notations of  the tropical Nevanlinna theory for generalized tropical meromorphic functions, and prove that any generalized tropical meromorphic function can be expressed as the tropical quotient of two generalized tropical entire functions (Proposition \ref{p11}). In the third section, we obtain the Poisson-Jensen formula for generalized tropical meromorphic function (Theorem \ref{T3.1}) and the first main theorem for generalized tropical meromorphic function with piecewise-linear target (Theorem \ref{T1.4}). In section 4, we prove a logarithmic derivative lemma for slope-admissible generalized tropical meromorphic functions (Lemma \ref{LDL}). In Section 5, we apply the theory to ultra-discrete equations, including the finite-order criterion for \eqref{eq13} and the existence of infinite-order solutions.\par

\section{Preliminaries} In this section, we will introduce the preliminaries of the tropical Nevanlinna theory (see references \cite{Halburd09, Risto15, Speyer09}).

First, we recall the definition of discontinuity points \cite{1}. Discontinuities are classified into three types: removable discontinuities, jump discontinuities and essential discontinuities. At discontinuity points of $f$, if the one-sided limits
\[
f(x_0^+) \in \mathbb{R},\qquad f(x_0^-) \in \mathbb{R}
\]
both exist and satisfy $f(x_0^+) \neq f(x_0^-)$, then $x_0$ is called a \textit{jump discontinuity}; If the above one-sided limits exist with $f(x_0^+) = f(x_0^-)$, but $f$ is undefined at $x_0$ or $f(x_0)\neq f(x_0^+)=f(x_0^-)$, then $x_0$ is called a \textit{removable discontinuity}; If one or both of the above one-sided limits fail to exist in $\mathbb{R}$, then $x_0$ is called an \textit{essential discontinuity}. Jump discontinuities and removable discontinuities are collectively referred to as \textit{first-kind discontinuities}, while  all other discontinuities (points where at least one one-sided limit fails to exist) are classified as \textit{second-kind discontinuities}. \par

\begin{lemma}\label{L2.2} All discontinuities of a generalized tropical meromorphic function are first-kind discontinuities.
\end{lemma}
\begin{proof} Let $f:\mathbb{R}\to\mathbb{R}$ be a generalized tropical meromorphic function. If $f$ is discontinuous at $x_0$, it follows from Definition \ref{D2.1} that there exist real numbers \(a,b\) with \(a<x_0<b\) such that $f$ is linear on \((a, x_0)\) and \((x_0, b)\). Hence, both $f(x_0^+)$ and $f(x_0^-)$ must exist. Therefore, every point of discontinuity of $f$ is necessarily a first-kind discontinuity.
\end{proof}

For a generalized tropical meromorphic function $f$, we define roots and poles in the same way as for tropical meromorphic function \cite{Risto15}. At the non-differentiable points of  \(f\), we say that $x$ is a \textit{pole} if
\[
\omega_f(x):=\lim_{\varepsilon\to 0_{+}}\left(f'(x + \varepsilon)-f'(x - \varepsilon)\right)<0
\] and a \textit{root} if
\(
\omega_f(x)>0.
\) The multiplicity of a root or a pole of $f$ is \(\tau_f(x)=|\omega_f(x)|\).\par

\begin{definition}[Generalized tropical entire function]\label{D2.3} Let $f: \mathbb{R}\to\mathbb{R}$ a generalized tropical meromorphic function. The function $f$ is said to be generalized tropical entire if it has no poles, right-negative jump points, or left-negative jump points.
\end{definition}

Here, for a generalized tropical meromorphic function \(f\), denote the \textit{tropical proximity function} by
\[
m(r,f):=\frac{1}{2}\sum_{\sigma=\pm1}f^{+}(\sigma r)=\frac{f^{+}(r)+f^{+}(-r)}{2},
\] where \(f^{+}(x)=\max\{f(x),0\}\), and the \textit{tropical counting function} is defined as
\[
N(r,f):=\frac{1}{2}\int_{0}^{r}n(t,f)dt=\frac{1}{2}\sum_{|b_v|<r}|\omega_f(b_v)|(r - |b_{\nu}|),
\] where
\[
n(t,f)=\sum_{\substack{|s|\leq t\\ \omega_f(s)<0}}|\omega_f(s)|.
\]

The generalized tropical meromorphic functions constructed in this section actually unify and contain both classical tropical meromorphic functions and the discontinuous piecewise linear functions considered by Halonen-Korhonen-Filipuk in \cite{Halonen23} as special cases, as shown below.

\begin{remark}\label{R1} (i). The piecewise linear functions in the sense of \cite{Halonen23} may have jump discontinuities, but at each discontinuity point they must be at least left-continuous or right-continuous. In this setting,  either $J^-(r, f)=0$ or $J^+(r, f)=0.$ Thus, in the one-sided continuous case, our theory reduces exactly to the results of \cite{Halonen23}.

(ii). If generalized tropical meromorphic function $f$ is continuous at every point, then it reduces to the classical tropical meromorphic function in which
 $$J^+(r, f)=J^-(r, f)=J(r, f)=0.$$
\end{remark}

Next, we prove the following result.

\begin{proposition}\label{L11} The class of generalized tropical meromorphic functions is closed under the tropical operations.
\end{proposition}
\begin{proof} For any two generalized tropical meromorphic functions $f$, $g$, by Definition \ref{D2.1}, we know that their partition points are either finite or countable, form a locally finite set, and $f$ and $g$ are linear on the open intervals between adjacent partition points. Clearly, the partition points of $f\oplus g$ and $f\otimes g$ are at most the union of the partition points of $f$ and $g$, so they are still finite or countable and form a locally finite set; on the open intervals between adjacent partition points, $f\oplus g$ and $f\otimes g$ are also linear. Therefore, $f\oplus g$ and $f\otimes g$ are generalized tropical meromorphic functions.
\end{proof}

By the closure under the tropical operations, it is straightforward to verify that the tropical proximity function, tropical counting function, right jump counting function and left jump counting function satisfy the following basic properties (see \cite[Lemma 3.2]{Risto15} for details).\par

\begin{lemma}\label{L2.3} Let $f$ and $g$ be two generalized tropical meromorphic functions. Then
    \[
    \begin{aligned}
    V(r,f^{\otimes \alpha}) &= \alpha V(r,f),~\alpha\in\mathbb{R},\\
    V(r,f \otimes g) &\leq V(r,f) + V(r,g), \\
    V(r,f \oplus g) &\leq V(r,f) + V(r,g)
    \end{aligned}
    \] hold for $V(r,f)=m(r,f)$, $V(r,f)=N(r,f)$, $V(r,f)=J^+(r,f)$ and $V(r,f)=J^-(r,f)$.
\end{lemma}

We recall the following measure estimate, which will be used in the proof of Lemma \ref{LDL} to estimate the exceptional set.
\begin{theorem}\cite[Theorem 3.3]{Halonen24}\label{T2.7}
Let $F(r)$ be a positive, nondecreasing and piecewise continuous function defined for $r_0\le r<\infty$, and let $\phi(r)$ be a positive, nondecreasing, continuous function defined for $r_0\le r<\infty$, and assume that $F(r)\ge e$ for $r\ge r_0$. Let $\xi(x)$ be a positive, nondecreasing, continuous function defined for $e\le x<\infty$. Let $C>1$ be a constant and let $E$ be the subset of $[r_0,\infty)$ defined by
\[
E=\left\{r\in [r_0,\infty): F\left(r+\frac{\phi(r)}{\xi(F(r))}\right)\ge C F(r)\right\}.
\]
Then for all $R<\infty$,
\[
\int_{E\cap [r_0,R]}\frac{dr}{\phi(r)}
\le \frac1{\xi(e)}+\frac1{\log C}\int_{e}^{F(R+)}\frac{dx}{x\xi(x)}.
\]
\end{theorem}

We can extend \cite[Proposition 3.3]{Risto16} analogously to the case of generalized tropical meromorphic functions, as given below.\par

\begin{proposition}\label{p11}
For any generalized tropical meromorphic function $f$, there exist generalized tropical entire functions $h$ and $g$ such that $f = h \oslash g$, where $h$ and $g$ have no common roots, right-positive jump points, or left-positive jump points.
\end{proposition}
\begin{proof} We will construct $h(x)$ and $g(x)$ separately to prove the assertion. Assume that
\[
f(x) = f_i(x) = a_ix+b_i,\quad x\in (x_{i-1},x_i),
\] where $\{x_k\}_{k\in\mathbb{Z}}\subset\mathbb{R}$ are all partition points of $f$ with the ordering $\cdots<x_{-2}<x_{-1}<x_0<x_1<x_2<\cdots$. We suppose
\begin{equation}\label{eq2.1}
g(x) = g_i(x) = c_ix+d_i,\quad x\in (x_{i-1},x_i),
\end{equation}
then
\begin{equation}\label{eq2.2}
h(x) = f(x)+g(x) = (a_i+c_i) x+(b_i+d_i),\quad x\in (x_{i-1},x_i).
\end{equation}
Now, we only need to determine the coefficients $c_i, d_i$ for $i\in\mathbb{Z}/\{0\}$. In fact, all of these coefficients are determined by $c_0$ and $d_0$, which can be chosen arbitrarily.

Given arbitrary $c_0, d_0\in \mathbb{R}$. Let
\begin{equation}\label{eq2.3}
\begin{split}
\omega_g(x_i)&:= (-\omega_f)^+(x_i),\\
sgn(x_i^+)\Omega_g(x_i^+)&:= \bigl(-sgn(x_i^+)\Omega_f(x_i^+)\bigl)^+,\\
sgn(x_i^-)\Omega_g(x_i^-)&:= \bigl(-sgn(x_i^-)\Omega_f(x_i^-)\bigl)^+,
\end{split}
\end{equation}
then
\begin{equation*}
\begin{split}
\omega_h(x_i) &= \omega_f(x_i)+\omega_g(x_i)=(\omega_f)^+(x_i),\\
sgn(x_i^+)\Omega_h(x_i^+)&= sgn(x_i^+)\bigl(\Omega_f(x_i^+)+\Omega_g(x_i^+)\bigl)=\bigl(sgn(x_i^+)\Omega_f(x_i^+)\bigl)^+,\\
sgn(x_i^-)\Omega_h(x_i^-)&= sgn(x_i^-)\bigl(\Omega_f(x_i^-)+\Omega_g(x_i^-)\bigl)=\bigl(sgn(x_i^-)\Omega_f(x_i^-)\bigl)^+,
\end{split}
\end{equation*} where $i\in\mathbb{Z}$. Obviously, the functions $h$ and $g$ have no poles, nor do they have right-negative jump points or left-negative jump points. Furthermore, $x_i$ cannot be a common root of $h$ and $g$, nor can it be a common right-positive jump point or left-positive jump point of $h$ and $g$.\par

It follows from \eqref{eq2.1} that
\begin{equation*}
\begin{split}
\omega_g(x_i)&=c_{i+1}-c_i,\\
c_ix_i+d_i+\Omega_g(x_i^+)+\Omega_g(x_i^-)&=c_{i+1}x_i+d_{i+1},\\
g(x_i)&=c_ix_i+d_i+\Omega_g(x_i^-),
\end{split}
\end{equation*} where $i\in\mathbb{Z}$. One can observe that, for $i\in\mathbb{Z}$,
\[
d_{i+1}-d_i=c_ix_i-c_{i+1}x_i+\Omega_g(x_i^+)+\Omega_g(x_i^-)=-\omega_g(x_i)x_i+\Omega_g(x_i^+)+\Omega_g(x_i^-).
\]
Hence, for $i\in\mathbb{N}$, we have
\begin{equation}\label{eq2.4}
\begin{split}
c_{i+1}&=c_0+\sum_{k=0}^{i}\omega_g(x_{k}),\\
d_{i+1}&=\sum_{k=0}^i(d_{k+1}-d_k)+d_0\\
&=\sum_{k=0}^i-\omega_g(x_k)x_k+\sum_{k=0}^i\bigl(\Omega_g(x_k^+)+\Omega_g(x_k^-)\bigl)+d_0.
\end{split}
\end{equation} Similarly, for $i\in\mathbb{Z}/\mathbb{N}$, we have
\begin{equation}\label{eq2.5}
\begin{split}
c_i&=c_0-\sum_{k=i}^{-1}\omega_g(x_{k}),\\
d_i&=\sum_{k=i}^{-1}(d_k-d_{k+1})+d_0\\
&=\sum_{k=i}^{-1}\omega_g(x_k)x_k-\sum_{k=i}^{-1}\bigl(\Omega_g(x_k^+)+\Omega_g(x_k^-)\bigl)+d_0.
\end{split}
\end{equation} Combining the equations \eqref{eq2.1}, \eqref{eq2.2}, \eqref{eq2.3}, \eqref{eq2.4} and \eqref{eq2.5}, we can see that $g(x)$ and $h(x)$ can be determined by $g_0$. Thus, the proof is complete.
\end{proof}

\section{Poisson-Jensen formula}
In this section, we prove the Poisson-Jensen formula for generalized tropical meromorphic functions.

\begin{theorem}[Poisson-Jensen formula]\label{T3.1} Suppose $f$ is a generalized tropical meromorphic function on $[-r,r]$ for some $r>0$. For any $x\in(-r,r)$, we let $(c_j)_{j=-p}^q$ be the increasing sequence of points in $(-r,r)$ such that $\omega_{f}\left(c_{j}\right)\neq 0$ for all $j\neq 0$, with $c_0=x$, and let $(k_j)_{j=1}^s$ be the increasing sequence of discontinuities of $f$ in $[-r,r]$. Then we have the Poisson-Jensen formula
\begin{equation*}
\begin{split}
f(x) =& \frac{1}{2}\bigl(f(r)+f(-r)\bigr) + \frac{x}{2r}\bigl(f(r)-f(-r)\bigr) \\
& +\frac{1}{2r}\sum_{j=-p}^{q}-\omega_{f}\left(c_{j}\right)\left(r^{2}-\left|c_{j}-x\right| r-c_{j} x\right) \\
& +\frac{x}{2r} \Bigg( \sum_{ -r\leq k_j^+\leq r}\kern-6pt-\Omega_{f}\left(k_{j}^+\right)+\sum_{-r\leq k_j^-\leq r}\kern-6pt-\Omega_{f}\left(k_{j}^-\right)\Bigg)\\
&+\frac12\Bigg( \sum_{-r\leq k_j^+\leq x}\kern-6pt\Omega_{f}(k_{j}^+) +\sum_{-r\leq k_j^-\leq x}\kern-6pt\Omega_{f}(k_{j}^-)-\sum_{ x\leq k_j^+\leq r}\kern-6pt\Omega_{f}(k_{j}^+)-\sum_{ x\leq k_j^-\leq r}\kern-6pt\Omega_{f}(k_{j}^-) \Bigg).
\end{split}
\end{equation*}
In particular, when $x=0$, we have the Jensen formula
\begin{align}\label{Jf}
f(0) = &\frac{1}{2}\bigl(f(r)+f(-r)\bigr) +\frac{1}{2}\sum_{j=-p}^{q}-\omega_{f}\left(c_{j}\right)\left(r-\left|c_{j}\right|\right) \notag\\
&+\frac12\Bigg( \sum_{-r\leq k_j^+\leq 0}\kern-6pt\Omega_{f}(k_{j}^+) +\sum_{-r\leq k_j^-\leq 0}\kern-6pt\Omega_{f}(k_{j}^-)-\sum_{ 0\leq k_j^+\leq r}\kern-6pt\Omega_{f}(k_{j}^+)-\sum_{ 0\leq k_j^-\leq r}\kern-6pt\Omega_{f}(k_{j}^-) \Bigg).
\end{align}
\end{theorem}
\begin{figure}[H]
  \centering
  \includegraphics[width=0.6\linewidth]{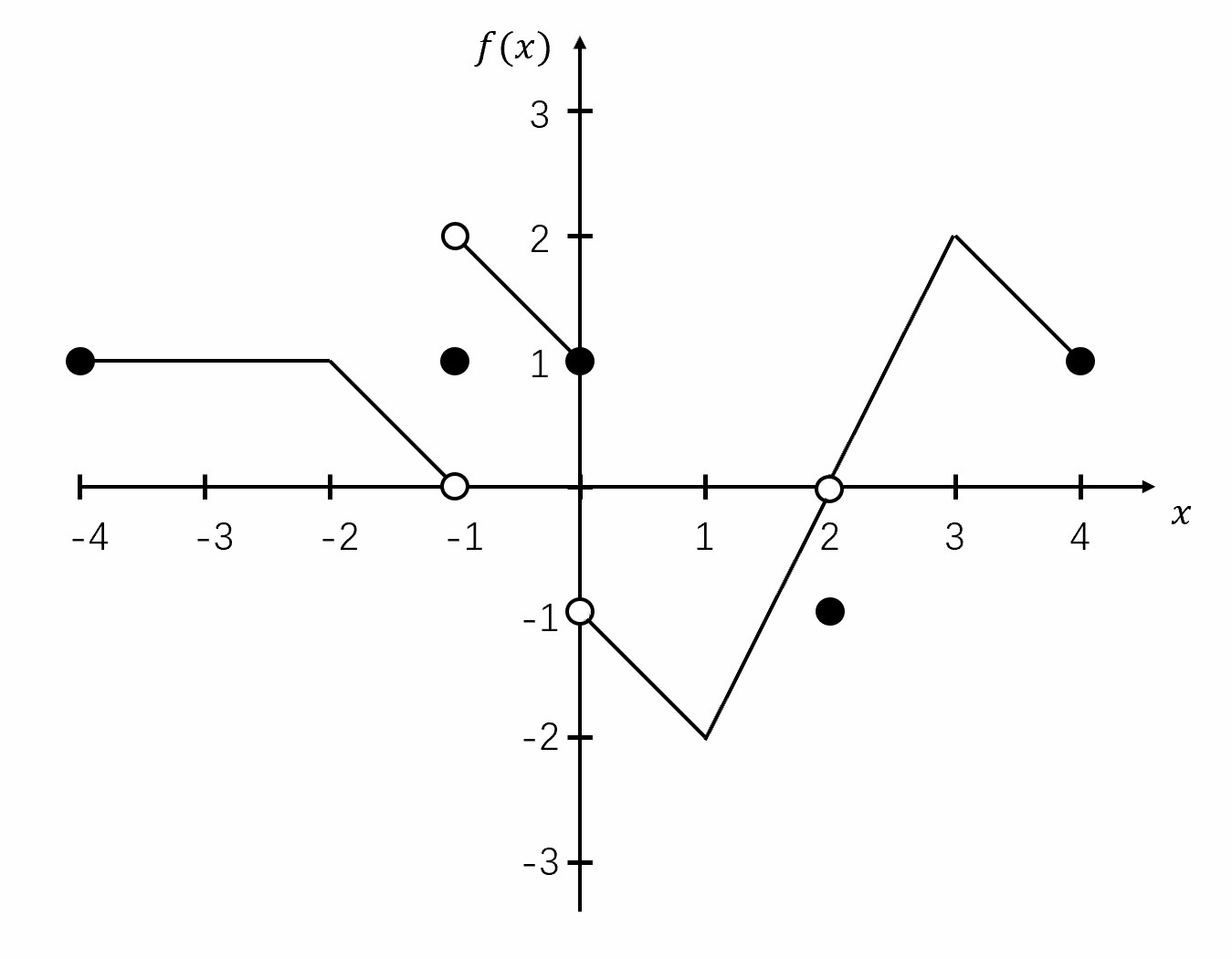}\\
  \caption{Function $f$ in Example 3.2.}\label{111.jpg}
\end{figure}
Next, we provide an example to verify this result.
\begin{example}\label{L3.2} Let $r=4$ and define $f:[-4,4]\to\mathbb{R}$,
\[
f(x)=
\begin{cases}
1, & -4\leq x < -2, \\
-x-1, & -2\leq x< -1, \\
1, & x=-1,\\
-x+1, & -1< x\leq 0,\\
-x-1, & 0< x< 1,\\
2x-4, & 1\leq x < 2,\\
-1, & x=2\\
2x-4, & 2< x \leq 3,\\
-x+5, & 3< x\leq 4.
\end{cases}
\]

This function has poles at $x=-2$ and $x=3,$   of multiplicities $1$ and $3$ respectively, and a root of multiplicity $3$ at $x=1.$ In addition, the function $f$ has three discontinuities: $x=-1,$ $x=0$ and $x=2$. Specifically, $x=-1$ is a step, this point is both a right-negative jump with height 1 and a left-negative jump with height 1;  $x=0$ is a jump whose assigned value lies outside the interval between the one-sided limits, this point is a right-negative jump with height 2; $x=2$ is a removable discontinuity, this point is a right-positive jump with height 1 and a left-negative jump with height 1.\par

Let $x\in(-1,0)$, by Poisson-Jensen formula, then we have
\begin{align*}
&f(x)=\frac12\bigl(1+1\bigr)+\frac{x}{8}\bigl(1-1\bigr)\\
& +\frac{1}{8}\left(16-4\left|-2-x\right|+2 x\right)-\frac{1}{8}3\left(16-4\left|1-x\right|- x\right)+\frac{1}{8}3\left(16-4\left|3-x\right|- 3x\right)   \\
& +\frac{x}{8}\bigl(-(1+(-2)+1)-(1+(-1))\bigl)+\frac12\bigl((1+1)-(-2+1+(-1))\bigl)\\
&\quad\quad=-x+1.
\end{align*}
\end{example}
\begin{figure}[H]
  \centering
  \includegraphics[width=1\linewidth]{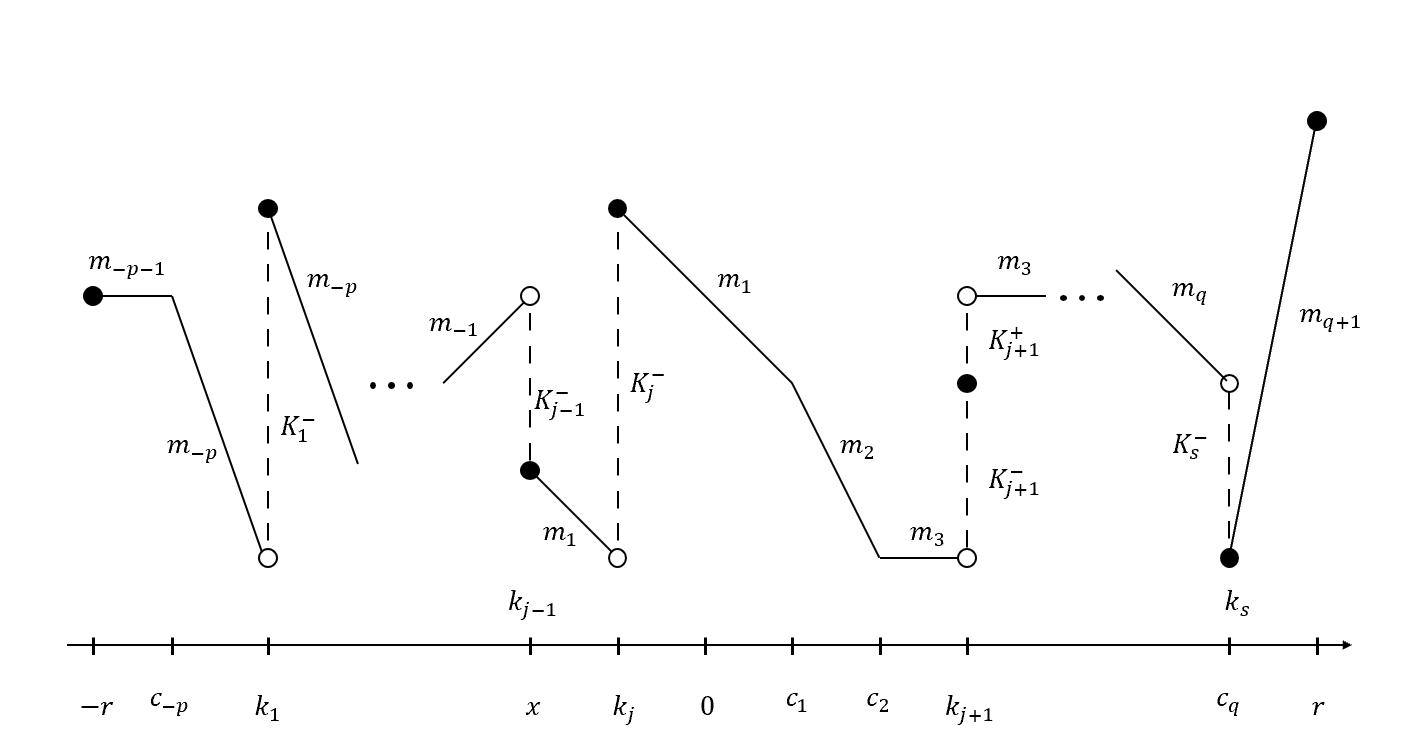}\\
  \caption{Notation in the proof of the Poisson-Jensen formula.}\label{11.jpg}
\end{figure}

\noindent{\textit{Proof of Theorem \ref{T3.1}.}} Denote
\[
m_{j-1} := f'(c_j^-)
\]
for $j=-p,\ldots,0$ and
\[
m_{j+1} := f'(c_j^+)
\]
for $j=0,\ldots,q$. Define
\[
K_j^+ = f(k_j^+) - f(k_j),\qquad K_j^- = f(k_j)-f(k_j^-)
\]
for $j=1,\ldots,s$. Also define the index sets
\[
J^+ = \{j\in\{1,\ldots,s\} : (x\le k_j\le r)\}
\]
and
\[
J^- = \{j\in\{1,\ldots,s\} : (-r\le k_j\le x)\}.
\]
A geometric observation implies that
\begin{align*}
f(r) - f(x) =& m_1(c_1 - x) + m_2(c_2 - c_1) + \cdots + m_q(c_q - c_{q-1}) + m_{q+1}(r - c_q) \\
&+ \sum_{\substack{j\in J^+\\ x\leq k_j^+\leq r}}\kern-6pt K_j^+ +\sum_{\substack{j\in J^+\\ x\leq k_j^-\leq r}}\kern-6pt K_j^-  \\
=& -m_1x + (m_1 - m_2)c_1 + \cdots + (m_q-m_{q+1})c_q + m_{q+1}r  \\
&+ \sum_{ x\leq k_j^+\leq r}\kern-6pt \bigl(f(k_j^+) - f(k_j)\bigl) +\sum_{ x\leq k_j^-\leq r}\kern-6pt \bigl(f(k_j)-f(k_j^-)\bigl) \\
=& m_1(r - x) - \sum_{j=1}^q (m_j - m_{j+1})(r - c_j) \\
&+ \sum_{x\leq k_j^+\leq r}\kern-6pt \bigl(f(k_j^+) - f(k_j)\bigl) +\sum_{x\leq k_j^-\leq r}\kern-6pt \bigl(f(k_j)-f(k_j^-)\bigl).
\end{align*}
Similarly,
\begin{align*}
f(x) - f(-r) =& m_{-1}(r + x) + \sum_{j=1}^p (m_{-j-1} - m_{-j})(r + c_{-j}) \\
&+ \sum_{ -r\leq k_j^+\leq x}\kern-6pt \bigl(f(k_j^+) - f(k_j)\bigl) +\sum_{-r\leq k_j^-\leq x}\kern-6pt \bigl(f(k_j)-f(k_j^-)\bigl).
\end{align*}
Multiplying the above two equalities by $(r+x)$ and $(r-x)$, respectively, and subtracting, we obtain
\[
\begin{aligned}
2 r f(x)
= &r(f(r)+f(-r))+x(f(r)-f(-r))+(m_{-1}-m_1)(r^2-x^2) \\
& +\sum_{j=1}^{p}\left(m_{-j-1}-m_{-j}\right)\left(r^{2}-\left(x-c_{-j}\right) r-c_{-j} x\right) \\
& +\sum_{j=1}^{q}\left(m_{j}-m_{j+1}\right)\left(r^{2}-\left(c_{j}-x\right) r-c_{j} x\right) \\
& +\Bigg( \sum_{ -r\leq k_j^+\leq x}\kern-6pt \bigl(f(k_j^+) - f(k_j)\bigl) +\sum_{-r\leq k_j^-\leq x}\kern-6pt \bigl(f(k_j)-f(k_j^-)\bigl)\Bigg)(r-x)\\
& -\Bigg( \sum_{x\leq k_j^+\leq r}\kern-6pt \bigl(f(k_j^+) - f(k_j)\bigl) +\sum_{x\leq k_j^-\leq r}\kern-6pt \bigl(f(k_j)-f(k_j^-)\bigl)\Bigg)(r+x)\\
=& r(f(r)+f(-r))+x(f(r)-f(-r)) \\
& +\sum_{j=-p}^{q}-\omega_{f}\left(c_{j}\right)\left(r^{2}-\left|c_{j}-x\right| r-c_{j} x\right) \\
& +x \Bigg( \sum_{ -r\leq k_j^+\leq r}\kern-6pt-\Omega_{f}\left(k_{j}^+\right)+\sum_{-r\leq k_j^-\leq r}\kern-6pt-\Omega_{f}\left(k_{j}^-\right)\Bigg)\\
&+r\Bigg( \sum_{-r\leq k_j^+\leq x}\kern-6pt\Omega_{f}(k_{j}^+) +\sum_{-r\leq k_j^-\leq x}\kern-6pt\Omega_{f}(k_{j}^-)-\sum_{ x\leq k_j^+\leq r}\kern-6pt\Omega_{f}(k_{j}^+)-\sum_{ x\leq k_j^-\leq r}\kern-6pt\Omega_{f}(k_{j}^-) \Bigg).
\end{aligned}
\]
Thus, the proof is complete.\qed

Combining the definitions of the tropical counting function and the total jump counting function, we can change Jensen formula \eqref{Jf} as
\begin{theorem}[Jensen formula]\label{T3.2} Let $f$ be a generalized tropical meromorphic function in $\mathbb{R}$. Then
\begin{equation}\label{eq11}
\begin{split}
f(0)=\frac{1}{2}[f(r)+f(-r)]-N(r, 1_{\mathbb{T}}\oslash f) + N(r, f)-J(r, 1_{\mathbb{T}}\oslash f) + J(r, f)
\end{split}
\end{equation} or equivalently,
\[
T(r, f)=T(r, 1_{\mathbb{T}}\oslash f)+f(0)
\] for any $r>0$.
\end{theorem}
\begin{proof} Here, we denote the distinct roots and poles of $f$ in $(-r,r)$ by $a_\mu$ and $b_\nu$, respectively; the distinct right-positive and right-negative jumps of $f$ in $[-r,r]$ by $\alpha_\mu$ and $\beta_\nu$, respectively; and the distinct left-positive and left-negative jumps of $f$ in $[-r,r]$ by $\xi_s$ and $\zeta_t$, respectively. Then we have
\begin{align}\label{eq888}
\frac12\sum_{j=-p}^{q}-\omega_{f}\left(c_{j}\right)\left(r-\left|c_{j}\right|\right) &=- \frac12\sum_{|a_\mu|<r} \tau_f(a_\mu)(r-|a_\mu|) +  \frac12\sum_{|b_\nu|<r} \tau_f(b_\nu)(r-|b_\nu|)\notag\\
&=-N(r, 1_{\mathbb{T}}\oslash f) + N(r, f)
\end{align} and we can see that
\begin{align}\label{eq999}
&\sum_{-r\leq k_j^+\leq 0}\kern-6pt\Omega_{f}(k_{j}^+) +\sum_{-r\leq k_j^-\leq 0}\kern-6pt\Omega_{f}(k_{j}^-)-\sum_{ 0\leq k_j^+\leq r}\kern-6pt\Omega_{f}(k_{j}^+)-\sum_{ 0\leq k_j^-\leq r}\kern-6pt\Omega_{f}(k_{j}^-)\notag\\
= &-\kern-6pt\sum_{-r \leq \alpha_{\mu}^+ \leq 0}\kern-6pt h_{f}\left(\alpha_{\mu}^+\right)+\sum_{-r \leq \beta_{\nu}^+ \leq 0}\kern-6pt h_{f}\left(\beta_{\nu}^+\right) -\Bigg(\sum_{0 \leq \alpha_{\mu}^+ \leq r}\kern-6pt h_{f}\left(\alpha_{\mu}^+\right)-\sum_{0 \leq \beta_{\nu}^+ \leq r}\kern-6pt h_{f}\left(\beta_{\nu}^+\right)\Bigg) \notag\\
& -\kern-6pt\sum_{-r \leq \xi_{s}^- \leq 0}\kern-6pt h_{f}\left(\xi_{s}^-\right)+\sum_{-r \leq \zeta_{t}^- \leq 0}\kern-6pt h_{f}\left(\zeta_{t}^-\right)-\Bigg(\sum_{0 \leq \xi_{s}^- \leq r}\kern-6pt h_{f}\left(\xi_{s}^-\right)-\sum_{0\leq \zeta_{t}^- \leq r}\kern-6pt h_{f}\left(\zeta_{t}^-\right)\Bigg) \\
=&- \sum_{|\alpha_\mu^+|\le r} h_f(\alpha_\mu^+) + \sum_{|\beta_\nu^+|\le r} h_f(\beta_\nu^+) -\sum_{|\xi_s^-|\le r} h_f(\xi_s^-) +  \sum_{|\zeta_t^-|\le r} h_f(\zeta_t^-)\notag\\
=&2\bigl(-J^+(r, 1_{\mathbb{T}}\oslash f) + J^+(r, f)-J^-(r, 1_{\mathbb{T}}\oslash f) + J^-(r, f)\bigl)\notag\\
=&2\bigl(-J(r, 1_{\mathbb{T}}\oslash f) + J(r, f)\bigl)\notag.
\end{align} Combining \eqref{Jf}, \eqref{eq888}, and \eqref{eq999}, the proof is complete.
\end{proof}

Next, we give the proof of the first main theorem (Theorem \ref{T1.4}) using the argument in \cite[Theorem 2.8]{Halonen23}.\vskip 6pt

\noindent{\bf Proof of Theorem \ref{T1.4}.} By the Jensen formula \eqref{eq11}, we have
\begin{align*}
T&\left(r,\frac{1_{\mathbb{T}}}{f\oplus a}\oslash\right)
-\bigl(T(r,f)-N(r,f)+N(r,f\oplus a)-J(r,f)+J(r,f\oplus a)\bigl) \\
&= T(r,f\oplus a)-\bigl(m(r,f)+N(r,f\oplus a)+J(r,f\oplus a)\bigr)-f(0)\oplus a(0) \\
&= m(r,f\oplus a)-m(r,f)-f(0)\oplus a(0).
\end{align*}
By the properties of the proximity function we obtain
\[
m(r,f\oplus a)-m(r,f)\le m(r,f)+m(r,a)-m(r,f)=m(r,a)
\]
and
\[
m(r,f\oplus a)-m(r,f)\ge m(r,f)-m(r,f)=0.
\]\qed

If $m(r,a)=o(T(r,f))$, where $r$ approaches infinity, then \eqref{eq7} can be written in the form
\begin{align*}\label{eq77}
T\left(r,\frac{1_{\mathbb{T}}}{f\oplus a}\oslash\right)= T(r,f)-N(r,f)+N(r,f\oplus a)-J(r,f)+J(r,f\oplus a)+o(T(r,f)),
\end{align*}
where $r$ approaches infinity.

\section{Logarithmic derivative lemma} According to the definition of the tropical characteristic function for a generalized tropical meromorphic function $f$, we can see that $T(r,f)$ is piecewise-linear in $[0, \infty)$. However, $T(r,f)$ is not necessarily non-decreasing or convex. We first give the following example to illustrate that $T(r,f)$ is not necessarily non-decreasing.

\begin{figure}[H]
  \centering
  \includegraphics[width=0.7\linewidth]{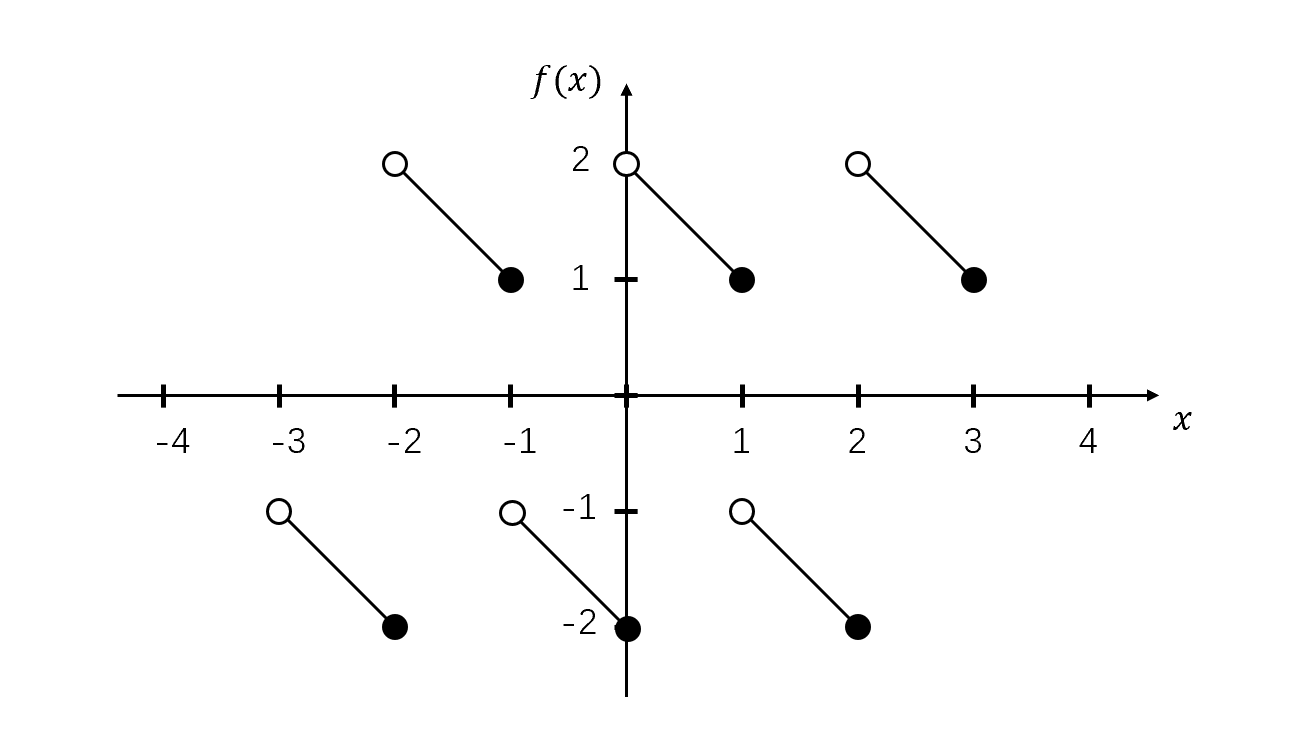}\\
  \caption{Function $f$ in Example \ref{E4.1}.}\label{2.jpg}
\end{figure}

\begin{example}\label{E4.1} For every $n\in\mathbb{Z}$, take a generalized tropical meromorphic function
\[
f(x)=
\begin{cases}
-x+2n+2, & 2n < x \leq 2n+1, \\
-x+2n, & 2n+1 < x \leq 2n+2.
\end{cases}.
\] From Figure \ref{2.jpg}, it is clear that the function $f$ is discontinuous at $x=n,$ $n\in\mathbb{Z}$ and $N(r,f)\equiv 0$ for all $r$. By a simple calculation, we can obtain
\begin{align*}
J(r,f)&=J^+(r,f)+J^-(r,f)=3n\\
m(r,f)&=\frac{-r+2n+2}{2}
\end{align*} for $2n< r< 2n+1$, where $ n\in \mathbb{N}$. Therefore,
\[
T(r,f)=-\frac{r}{2}+4n+1, \quad 2n< r< 2n+1 \quad n\in \mathbb{N}.
\]  Clearly, it is decreasing on each interval $(2n, 2n+1),$ where $ n\in \mathbb{N}$.
\end{example}

Next, we prove that $T(r, f)$ is nondecreasing when $f$ is slope-admissible.

\begin{theorem}\label{T4.3} Assume that generalized tropical meromorphic function $f$ is slope-admissible, then there exists a constant $r_0\in\mathbb{R}$, such that the tropical characteristic function $T(r, f)$ of $f$ is non-decreasing function for any $r>r_0$.
\end{theorem}
\begin{proof} Define $r_0$ as the minimum value of $\xi$ such that there exist \(x_1,x_2\in[-\xi, \xi]\) satisfying \(f'(x_1^+)\geq 0\) and \(f'(x_2^+)\leq 0\). To prove the monotonicity of $T(r, f)$ in $[r_0, \infty)$, we consider the following four cases.\par

Case 1: If $f(-r) < 0$, $f(r) < 0$, then $m(r, f) = 0$. Clearly, $m(r, f)$ is non-decreasing as $r$ increases from $r$ to $r + \varepsilon$. Since $N(r, f)$, $J^+(r,f)$ and $J^-(r,f)$ are increasing, it follows that $T(r, f)$ is also non-decreasing as $r$ increases from $r$ to $r + \varepsilon$.\par

Case 2: If $f(-r) \geq 0$, $f(r) \geq 0$, then $-f(-r) \leq 0$, $-f(r) \leq 0$, which implies $m(r, -f) = 0$. Clearly, $m(r, -f)$ is non-decreasing as $r$ increases from $r$ to $r + \varepsilon$, hence $T(r, -f)$ is non-decreasing. By the Jensen formula
\[
T(r, f) = T(r, -f) + f(0).
\] Hence, $T(r, f)$ is also non-decreasing.\par

Case 3: If $f(-r) < 0$ and $f(r) \geq 0$, then $m(r, f) = \frac{f(r)}{2}$. We know that \(f'(x_1^+)\geq 0\), where $x_1\in[-r_0, r_0]\subset [-r, r]$. Define an increasing sequence $(c_j)$, $j=1,\dots,q$ in $(x_1, r]$ in the following way. Let the points in this sequence be the points in $(x_1, r]$, at which the derivative of $f$ does not exist.\par

Next, we first assume that $r$ is a right continuous point (i.e., $f(r)=f(r^+)$). Then we have
\begin{eqnarray}\label{eq88} f'(r^+)-f'(x_1^+)&=&f'(r^+)-f'(c_q^+)+f'(c_q^+)-f'(c_{q-1}^+)+\cdots+f'(c_3^+)-f'(c_2^+)\nonumber\\
&~&+f'(c_2^+)-f'(c_1^+)+f'(c_1^+)-f'(x_1^+)\nonumber\\
&=&\sum_{j=1}^q\omega_f(c_j)\\
&=&\sum_{x_1 \leq a_v \leq r}\tau_f(a_v)-\sum_{x_1 \leq b_v \leq r}\tau_f(b_v)\nonumber.
\end{eqnarray}
Combining the above equation \eqref{eq88}, we can obtain
\begin{equation}\label{eq8}
\begin{split}
T'(r^+, f) &\geq \frac{1}{2}\bigl(f'(r^+) + n(r^+, f)  \bigr)\\
&\geq\frac{1}{2}\bigl(f'(r^+) + \sum_{x_1 \leq b_v \leq r}\tau_f(b_v)\bigr)\\
&=\frac{1}{2}\bigl(\sum_{x_1 \leq a_v \leq r}\tau_f(a_v) + f'(x_1^+)\bigr)\geq0,
\end{split}
\end{equation}
hence $T(r, f)$ is increasing at $r$.\par

Assume that $r$ is a right discontinuity point (i.e., $f(r)\neq f(r^+)$), then $f(r^+)=f(r)+\Omega_f(r^+)$. We discuss it in three subcases.\par

Subcase 3.1: If $f(r^+)\leq f(r)$, then $\Omega_f(r^+)<0$, this implies that $r$ is a right-negative jump of $f$, hence we have
\begin{equation*}
\begin{split}
m(r,f)+J^+(r,f)&=\frac{1}{2}\bigl(f(r)+\Omega_f(r^+)\bigl)-\frac{1}{2}\Omega_f(r^+)+J^+(r,f)\\
&=\frac12f(r^+)+\frac{1}{2}|\Omega_f(r^+)|+J^+(r,f)\\
&\leq m(r+\varepsilon ,f)+J^+(r+\varepsilon ,f)
\end{split}
\end{equation*} as $r$ increases from $r$ to $r + \varepsilon$. Since $N(r, f)$ and $J^-(r,f)$ are increasing, it follows that $T(r, f)$ is non-decreasing.\par

Subcase 3.2: If $f(r^+)>f(r)$, then $ m(r+\varepsilon ,f)>m(r,f)$. Hence $T(r, f)$ is non-decreasing as $r$ increases from $r$ to $r + \varepsilon$.\par

Case 4: Finally, if $f(-r) \geq 0$ and $f(r) < 0$, then $m(r, f) = \frac{f(-r)}{2}$. We know that \(f'(x_2^+)\leq 0\), where $x_2\in[-r_0, r_0]\subset [-r, r]$. Define an increasing sequence $(c_j)$, $j=1,\dots,p$ in $[-r, x_2)$ in the following way. Let the points in this sequence be the points in $[-r, x_2)$, at which the derivative of $f$ does not exist. Assume that $-r$ is a left continuous point (i.e., $f(-r)=f((-r)^-)$). Then we have
\begin{eqnarray}\label{eq99} f'(x_2^+)-f'((-r)^-)&=&f'(x_2^+)-f'(c_p^+)+f'(c_p^+)-f'(c_{p-1}^+)+\cdots+f'(c_3^+)-f'(c_2^+)\nonumber\\
&~&+f'(c_2^+)-f'(c_1^+)+f'(c_1^+)-f'((-r)^-)\nonumber\\
&=&\sum_{j=1}^p\omega_f(c_j)\\
&=&\sum_{-r \leq a_v \leq x_2}\tau_f(a_v)-\sum_{-r \leq a_v \leq x_2}\tau_f(b_v)\nonumber.
\end{eqnarray}
Combining the above equation \eqref{eq99}, we can obtain
\begin{equation*}
\begin{split}
T'(r^+, f) &\geq \frac{1}{2}\bigl(-f'((-r)^-) + n(r^+, f)  \bigr)\\
&\geq\frac{1}{2}\bigl(-f'((-r)^-) + \sum_{-r \leq a_v \leq x_2}\tau_f(b_v)\bigr)\\
&=\frac{1}{2}\bigl(\sum_{-r \leq a_v \leq x_2}\tau_f(a_v) - f'(x_2^+)\bigr)\geq0,
\end{split}
\end{equation*}
hence $T(r, f)$ is increasing at $r$.\par

Assume that $-r$ is a left discontinuity point (i.e., $f(-r)\neq f((-r)^-)$), then $f(-r)=f((-r)^-)+\Omega_f((-r)^-)$. By arguments similar to those in Subcases 3.1 and 3.2, we also deduce that $T(r, f)$ is non-decreasing as $r$ increases from $r$ to $r + \varepsilon$.\par

Thus, the proof is complete.
\end{proof}

Next, we show that even if Definition \ref{D4.2} is satisfied, that is, $f$ is slope-admissible generalized tropical meromorphic function, the tropical characteristic function $T(r,f)$ need not be convex due to the presence of jumps. The following example illustrates this.

\begin{figure}[H]
  \centering
  \includegraphics[width=0.7\linewidth]{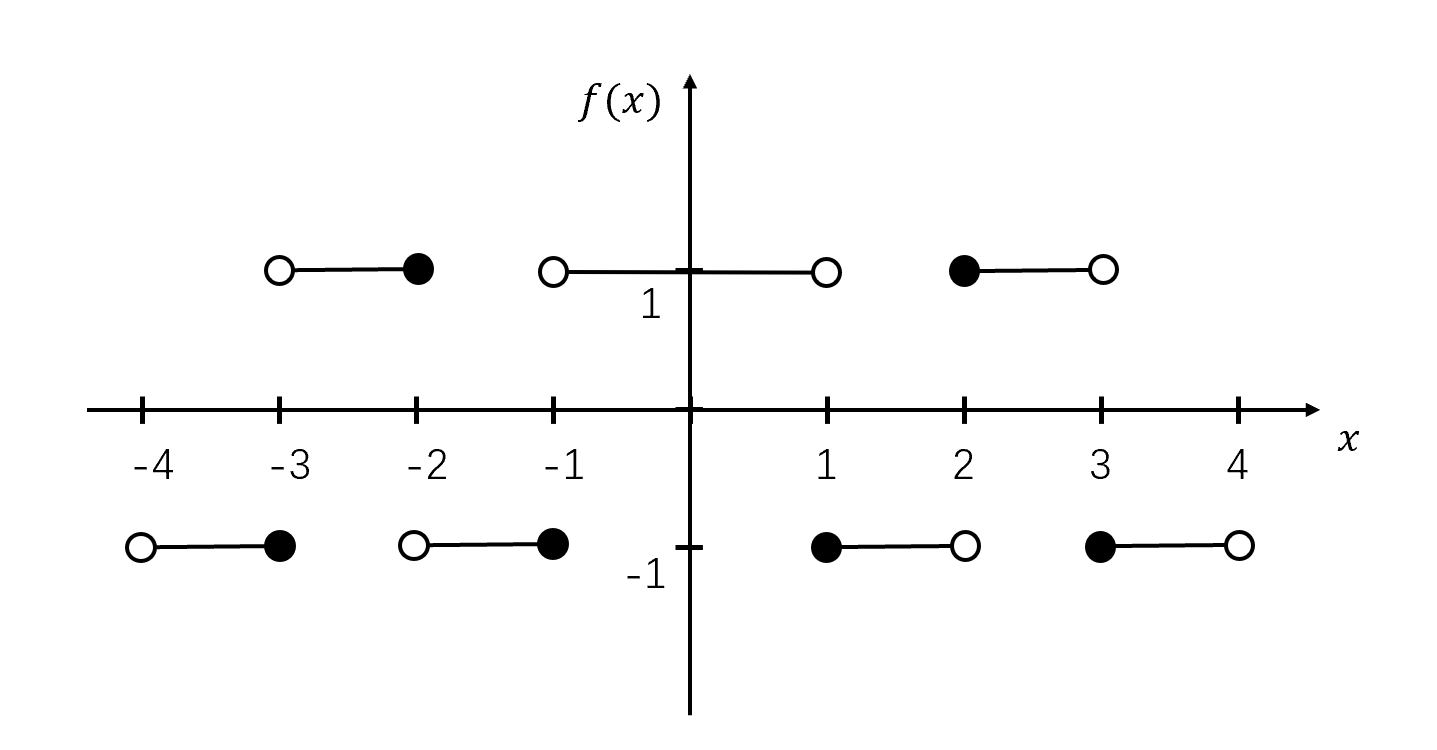}\\
  \caption{Functions in Example \ref{E4.4}.}\label{3.jpg}
\end{figure}

\begin{example}\label{E4.4} For every $n\in\mathbb{N}\cup \{0\}$, define \[
g(x)=
\begin{cases}
1, & 2n \leq x < 2n+1, \\
-1, & 2n+1 \leq x < 2n+2,
\end{cases} \text{ and } f(x)=
\begin{cases}
g(x), &  x\geq 0, \\
g(-x), & x<0,
\end{cases}.
\]  It follows that $$J(r,f)=J^+(r,f)+J^-(r,f)=2n$$ for $2n-1\leq r< 2n+1$, where $ n\in \mathbb{N}.$ Moreover, \[
m(r,f)=
\begin{cases}
1, & 2n \leq r < 2n+1, \\
0, & 2n+1 \leq r < 2n+2,
\end{cases}.
\] Clearly, $N(r,f)\equiv 0$ for all $r$. Therefore,
\[
T(r,f)=n+1, \quad n \leq r <n+1 \quad n\in \mathbb{N}\cup \{0\}.
\] Now, let $x=\frac12$ and $y=\frac32$. Then $$T(\frac{x+y}{2}, f)=2> \frac{T(x,f)+T(y,f)}{2}=\frac{3}{2}.$$ Hence, $T(r,f)$ is not convex.
\end{example}

The following theorem is a slight modification of \cite[Theorem 3.1]{Halonen24}, in which the domain $[0,\infty)$ is replaced by $[r_0,+\infty)$ for some constant $r_0>0$. The proof proceeds exactly as in \cite{Halonen24}, with only trivial adjustments for the starting point of sequences and integration intervals, so we omit the detailed proof here.

\begin{theorem}\label{T4.4}
Let $s>0$, $r_0>0$ and let $T\colon [r_0,\infty)\to [0,\infty)$ be a generalized tropical meromorphic function such that
\begin{equation}
T(r)\le \frac{R}{r}T(R)
\label{eq:condition-modified}
\end{equation}
for all $R\ge r\ge r_0$. If the hyper-order of $T$ is strictly less than one, then
\[
T(r+s)=T(r)+o(T(r))
\]
holds outside a set of finite logarithmic measure as $r$ tends to infinity.
\end{theorem}

Next, we give an estimate of the proximity function of $f(x+c)\oslash f(x)$.

\begin{theorem}\label{T4.33} Let $c\in\mathbb{R}\setminus\{0\}$. If $f$ is a generalized tropical meromorphic function. Then, for all $\alpha>1$ and all $r>\max\{1, |c|\}$, we have
\[
\begin{aligned}
m(r, f(x+c)\oslash f(x))&\leq\frac{16|c|}{r+|c|}\frac{1}{\alpha-1}\left(T(\alpha(r+|c|), f)+\frac{|f(0)|}{2}\right)\\
&\quad + J(r+|c|,f)-J(r-|c|,f)+J(r+|c|,-f)-J(r-|c|,-f).
\end{aligned}
\]
\end{theorem}
\begin{proof} Take $\rho:=\frac{1}{2}(\alpha+1)(r+|c|)>r+|c|$.\par
Case 1. Suppose $c>0,$ then by the Poisson-Jensen formula in Theorem \ref{T3.1}, we have
\[
\begin{aligned}
f(&x+c)-f(x)=\frac{c}{2\rho}(f(\rho)-f(-\rho))\\
&+\frac{1}{2\rho} \sum_{|c_j|<\rho}\omega_{f}\left(c_{j}\right)\bigl((|c_{j} - x - c| - |c_{j} - x|)\rho + c c_{j}\bigr) \\
& - \frac{c}{2\rho} \Bigg( \sum_{ -\rho\leq k_j^+\leq \rho}\kern-6pt\Omega_{f}\left(k_{j}^+\right)+\sum_{-\rho\leq k_j^-\leq \rho}\kern-6pt\Omega_{f}\left(k_{j}^-\right)\Bigg)\\
&+\frac12\Bigg( \sum_{x< k_j^+\leq x+c}\kern-6pt\Omega_{f}(k_{j}^+) +\sum_{x< k_j^-\leq x+c}\kern-6pt\Omega_{f}(k_{j}^-)+\sum_{ x\leq k_j^+ < x+c}\kern-6pt\Omega_{f}(k_{j}^+)+\sum_{ x\leq k_j^-< x+c}\kern-6pt\Omega_{f}(k_{j}^-) \Bigg).
\end{aligned}
\] Next, we further discuss the roots, poles, right-positive jumps, right-negative jumps, left-positive jumps and left-negative jumps of $f$, we can see that
\begin{align*}
\frac{1}{2\rho} \sum_{|c_j|<\rho}\omega_{f}\left(c_{j}\right)&\bigl((|c_{j} - x - c| - |c_{j} - x|)\rho + c c_{j}\bigr)\\
=&\frac{1}{2\rho} \sum_{|a_\mu|<\rho} \tau_f(a_\mu)\bigl((|a_{\mu} - x - c| - |a_{\mu} - x|)\rho + c a_{\mu}\bigr) \\
&- \frac{1}{2\rho} \sum_{|b_\nu|<\rho}\tau_f(b_\nu) \bigl((|b_{\nu} - x - c| - |b_{\nu} - x|)\rho + c b_{\nu}\bigr)
\end{align*} and
\begin{align*}
\sum_{ -\rho\leq k_j^+\leq \rho}\kern-6pt\Omega_{f}\left(k_{j}^+\right)+\sum_{-\rho\leq k_j^-\leq \rho}\kern-6pt\Omega_{f}\left(k_{j}^-\right)=&\sum_{|\alpha_{\mu}^+ |\leq \rho} \Omega_{f}\left(\alpha_{\mu}^+\right)+\sum_{|\beta_{\nu}^+ |\leq \rho}\Omega_{f}\left(\beta_{\nu}^+ \right)\\
&+\sum_{|\xi_{s}^- |\leq \rho} \Omega_{f}\left(\xi_{s}^-\right)+\sum_{|\zeta_{t}^- |\leq \rho}\Omega_{f}\left(\zeta_{t}^- \right).
\end{align*}
We can also see that
\[
\begin{aligned}
\sum_{x< k_j^+\leq x+c}\kern-6pt\Omega_{f}(k_{j}^+)+\sum_{ x\leq k_j^+ < x+c}\kern-6pt\Omega_{f}(k_{j}^+)=&\sum_{k_j^+\in[x,x+c]}\kern-6pt\Omega_{f}(k_{j}^+)+\sum_{k_j^+\in(x,x+c)}\kern-6pt\Omega_{f}(k_{j}^+)\\
=&\sum_{\alpha_{\mu}^+\in[x,x+c]}\kern-6pt\Omega_{f}\left(\alpha_{\mu}^+\right)+\sum_{\alpha_{\mu}^+\in(x,x+c)}\kern-6pt\Omega_{f}\left(\alpha_{\mu}^+\right)\\
&+\sum_{\beta_{\nu}^+\in[x,x+c]}\kern-6pt\Omega_{f}\left(\beta_{\nu}^+\right)+\sum_{\beta_{\nu}^+\in(x,x+c)}\kern-6pt\Omega_{f}\left(\beta_{\nu}^+\right).
\end{aligned}
\] Similarly, we have
\[
\begin{aligned}
\sum_{ x< k_j^-\leq x+c}\kern-6pt\Omega_{f}(k_{j}^-)+\sum_{ x\leq k_j^-< x+c}\kern-6pt\Omega_{f}(k_{j}^-)=&\sum_{\xi_{s}^-\in[x,x+c]}\kern-6pt\Omega_{f}\left(\xi_{s}^-\right)+\sum_{\xi_{s}^-\in(x,x+c)}\kern-6pt\Omega_{f}\left(\xi_{s}^-\right)\\
&+\sum_{\zeta_{t}^-\in[x,x+c]}\kern-6pt\Omega_{f}\left(\zeta_{t}^-\right)+\sum_{\zeta_{t}^-\in(x,x+c)}\kern-6pt\Omega_{f}\left(\zeta_{t}^-\right).
\end{aligned}
\]
By combining these equations, we have
\[
\begin{aligned}
f(&x+c)-f(x)=\frac{c}{2\rho}(f(\rho)-f(-\rho))\\
&+ \frac{1}{2\rho} \sum_{|a_\mu|<\rho} \tau_f(a_\mu)\bigl((|a_{\mu} - x - c| - |a_{\mu} - x|)\rho + c a_{\mu}\bigr) \\
&- \frac{1}{2\rho} \sum_{|b_\nu|<\rho}\tau_f(b_\nu) \bigl((|b_{\nu} - x - c| - |b_{\nu} - x|)\rho + c b_{\nu}\bigr)\\
&- \frac{c}{2\rho}\left(\sum_{|\alpha_{\mu}^+ |\leq \rho} \Omega_{f}\left(\alpha_{\mu}^+\right)+\sum_{|\beta_{\nu}^+ |\leq \rho}\Omega_{f}\left(\beta_{\nu}^+ \right)+\sum_{|\xi_{s}^- |\leq \rho} \Omega_{f}\left(\xi_{s}^-\right)+\sum_{|\zeta_{t}^- |\leq \rho}\Omega_{f}\left(\zeta_{t}^- \right)\right)\\
&+\frac{1}{2}\left(\sum_{\alpha_{\mu}^+\in[x,x+c]}\kern -8pt\Omega_{f}\left(\alpha_{\mu}^+\right)+\sum_{\alpha_{\mu}^+\in(x,x+c)}\kern -8pt\Omega_{f}\left(\alpha_{\mu}^+\right)+
\sum_{\beta_{\nu}^+\in[x,x+c]}\kern -8pt\Omega_{f}\left(\beta_{\nu}^+\right)+\sum_{\beta_{\nu}^+\in(x,x+c)}\kern -8pt\Omega_{f}\left(\beta_{\nu}^+\right)\right)\\
&+\frac{1}{2}\left(\sum_{\xi_{s}^-\in[x,x+c]}\kern -8pt\Omega_{f}\left(\xi_{s}^-\right)+\sum_{\xi_{s}^-\in(x,x+c)}\kern -8pt\Omega_{f}\left(\xi_{s}^-\right)+
\sum_{\zeta_{t}^-\in[x,x+c]}\kern -8pt\Omega_{f}\left(\zeta_{t}^-\right)+\sum_{\zeta_{t}^-\in(x,x+c)}\kern -8pt\Omega_{f}\left(\zeta_{t}^-\right)\right).
\end{aligned}
\] By taking the proximity function on both sides and by the proof of \cite[Theorem 3.23]{Risto15}, we deduce that
\[
\begin{aligned}
&m(r,f(x+c)\oslash f(x))\le |c|\biggl(\frac1{\rho}\bigl(m(\rho,f)+m(\rho,-f)\bigl)+n(\rho,f)+n(\rho,-f)\biggr) \\
&+ m\biggl(r,-\frac{c}{2\rho}\Bigl(\sum_{|\alpha_{\mu}^+ |\leq \rho} \Omega_{f}\left(\alpha_{\mu}^+\right)+\sum_{|\beta_{\nu}^+ |\leq \rho}\Omega_{f}\left(\beta_{\nu}^+ \right)+\sum_{|\xi_{s}^- |\leq \rho} \Omega_{f}\left(\xi_{s}^-\right)+\sum_{|\zeta_{t}^- |\leq \rho}\Omega_{f}\left(\zeta_{t}^- \right)\Bigr)\biggr)\\
&+ m\biggl(r,\frac{1}{2}\bigl(\sum_{\alpha_{\mu}^+\in[x,x+c]}\kern -8pt\Omega_{f}\left(\alpha_{\mu}^+\right)+\sum_{\alpha_{\mu}^+\in(x,x+c)}\kern -8pt\Omega_{f}\left(\alpha_{\mu}^+\right)+
\sum_{\beta_{\nu}^+\in[x,x+c]}\kern -8pt\Omega_{f}\left(\beta_{\nu}^+\right)+\sum_{\beta_{\nu}^+\in(x,x+c)}\kern -8pt\Omega_{f}\left(\beta_{\nu}^+\right)\bigl)\biggr)\\
&+ m\biggl(r,\frac{1}{2}\bigl(\sum_{\xi_{s}^-\in[x,x+c]}\kern -8pt\Omega_{f}\left(\xi_{s}^-\right)+\sum_{\xi_{s}^-\in(x,x+c)}\kern -8pt\Omega_{f}\left(\xi_{s}^-\right)+
\sum_{\zeta_{t}^-\in[x,x+c]}\kern -8pt\Omega_{f}\left(\zeta_{t}^-\right)+\sum_{\zeta_{t}^-\in(x,x+c)}\kern -8pt\Omega_{f}\left(\zeta_{t}^-\right)\bigl)\biggr).
\end{aligned}
\]
Since the quantity below inside the proximity function is constant with respect to $x$,
\[
\begin{aligned}
m\biggl(&r,-\frac{c}{2\rho}\Bigl(\sum_{|\alpha_{\mu}^+ |\leq \rho} \Omega_{f}\left(\alpha_{\mu}^+\right)+\sum_{|\beta_{\nu}^+ |\leq \rho}\Omega_{f}\left(\beta_{\nu}^+ \right)+\sum_{|\xi_{s}^- |\leq \rho} \Omega_{f}\left(\xi_{s}^-\right)+\sum_{|\zeta_{t}^- |\leq \rho}\Omega_{f}\left(\zeta_{t}^- \right)\Bigr)\biggr)\\
&=\biggl(-\frac{c}{2\rho}\Bigl(\sum_{|\alpha_{\mu}^+|\le\rho}\Omega_f(\alpha_{\mu}^+)+\sum_{|\beta_\nu^+|\le\rho}\Omega_f(\beta_\nu^+)+\sum_{|\xi_{s}^- |\leq \rho} \Omega_{f}\left(\xi_{s}^-\right)+\sum_{|\zeta_{t}^- |\leq \rho}\Omega_{f}\left(\zeta_{t}^- \right)\Bigr)\biggr)^+ \\
&\le \frac{|c|}{2\rho}\biggl(\sum_{|\alpha_{\mu}^+|\le\rho}|\Omega_f(\alpha_{\mu}^+)|+\sum_{|\beta_\nu^+|\le\rho}|\Omega_f(\beta_\nu^+)|+\sum_{|\xi_{s}^- |\leq \rho} |\Omega_{f}\left(\xi_{s}^-\right)|+\sum_{|\zeta_{t}^- |\leq \rho}|\Omega_{f}\left(\zeta_{t}^- \right)|\biggr) \\
&=\frac{|c|}{\rho}\bigl(J^+(\rho,f)+J^+(\rho,-f)+J^-(\rho,f)+J^-(\rho,-f)\bigr).
\end{aligned}
\]
By the definition of the proximity function, we have
\[
\begin{aligned}
&m\biggl(r,\frac{1}{2}\biggr(\sum_{\alpha_{\mu}^+\in[x,x+c]}\kern -8pt\Omega_{f}\left(\alpha_{\mu}^+\right)+\sum_{\alpha_{\mu}^+\in(x,x+c)}\kern -8pt\Omega_{f}\left(\alpha_{\mu}^+\right)\biggr)\biggr) \\
& \quad =\frac14\Bigl(\sum_{\alpha_{\mu}^+\in[r,r+c]}\kern-12pt \Omega_{f}\left(\alpha_{\mu}^+\right)+\sum_{\alpha_{\mu}^+\in(r,r+c)}\kern-12pt\Omega_{f}\left(\alpha_{\mu}^+ \right)\Bigr)^+
+\frac14\Bigl(\sum_{\alpha_{\mu}^+\in[-r,-r+c]}\kern-16pt\Omega_f(\alpha_{\mu}^+)+\sum_{\alpha_{\mu}^+\in(-r,-r+c)}\kern-16pt\Omega_f(\alpha_{\mu}^+)\Bigr)^+ \\
& \quad \le \frac12\sum_{\alpha_{\mu}^+\in[r,r+c]}\kern-12pt|\Omega_f(\alpha_{\mu}^+)|
+\frac12\sum_{\alpha_{\mu}^+\in[-r,-r+c]}\kern-16pt|\Omega_f(\alpha_{\mu}^+)|\\
& \quad \le \frac12\bigl(\sum_{\alpha_{\mu}^+\in[r-c,r+c]}\kern-16pt|\Omega_f(\alpha_{\mu}^+)|
+\sum_{\alpha_{\mu}^+\in[-r-c,-r+c]}\kern-18pt|\Omega_f(\alpha_{\mu}^+)|\bigl)\\
& \quad \leq J^+(r+c,f)-J^+(r-c,f).
\end{aligned}
\]
Similarly,  we have
\[
\begin{aligned}
m\biggl(r,\frac{1}{2}\biggr(\sum_{\beta_{\nu}^+\in[x,x+c]}\kern -8pt\Omega_{f}\left(\beta_{\nu}^+\right)+\sum_{\beta_{\nu}^+\in(x,x+c)}\kern -8pt\Omega_{f}\left(\beta_{\nu}^+\right)\biggr)\biggr)&\leq J^+(r+c,-f)-J^+(r-c,-f),\\
m\biggl(r,\frac{1}{2}\biggr(\sum_{\xi_{s}^-\in[x,x+c]}\kern -8pt\Omega_{f}\left(\xi_{s}^-\right)+\sum_{\xi_{s}^-\in(x,x+c)}\kern -8pt\Omega_{f}\left(\xi_{s}^-\right)\biggr)\biggr)&\leq J^-(r+c,f)-J^-(r-c,f),
\end{aligned}
\]
and
\[
\begin{aligned}
m\biggl(r,\frac{1}{2}\biggr(\sum_{\zeta_{t}^-\in[x,x+c]}\kern -8pt\Omega_{f}\left(\zeta_{t}^-\right)+\sum_{\zeta_{t}^-\in(x,x+c)}\kern -8pt\Omega_{f}\left(\zeta_{t}^-\right)\biggr)\biggr)\leq J^-(r+c,-f)-J^-(r-c,-f).
\end{aligned}
\]

Combining the above equations, we have
\[
\begin{aligned}
&m(r,f(x+c)\oslash f(x))\\
\le |c|&\frac1{\rho}\Bigl(m(\rho,f)+m(\rho,-f)+J^+(\rho,f)+J^+(\rho,-f)+J^-(\rho,f)+J^-(\rho,-f)\Bigr)\\
+&|c|\bigl(n(\rho,f)+n(\rho,-f)\bigl) + J^+(r+c,f)-J^+(r-c,f)+J^+(r+c,-f)-J^+(r-c,-f)\\
+&J^-(r+c,f)-J^-(r-c,f)+J^-(r+c,-f)-J^-(r-c,-f)\\
\le |c|&\frac1{\rho}\Bigl(m(\rho,f)+m(\rho,-f)+J(\rho,f)+J(\rho,-f)\Bigr)+|c|\bigl(n(\rho,f)+n(\rho,-f)\bigl) \\
+& J(r+c,f)-J(r-c,f)+J(r+c,-f)-J(r-c,-f).
\end{aligned}
\]
Now by \cite[Lemma 3.21]{Risto15},
\[
n(\rho,f)+n(\rho,-f)\le \frac4{\alpha-1}\frac1{r+|c|}\Bigl(T(\alpha(r+|c|),f)+T(\alpha(r+|c|),-f)\Bigr).
\]
Since $1/(\alpha+1)\le 2/(\alpha-1)$, the Jensen formula implies that
\[
\begin{aligned}
m(r,&f(x+c)\oslash f(x))\\
&\le |c|\biggl(\frac2{(\alpha+1)(r+|c|)}\Bigl(T(\alpha(r+|c|),f)+T(\alpha(r+|c|),-f)\Bigr)\biggr) \\
&\quad + |c|\biggl(\frac4{\alpha-1}\frac1{r+|c|}\Bigl(T(\alpha(r+|c|),f)+T(\alpha(r+|c|),-f)\Bigr)\biggr) \\
&\quad + J(r+c,f)-J(r-c,f)+J(r+c,-f)-J(r-c,-f)\\
&\leq\frac{16|c|}{r+|c|}\frac{1}{\alpha-1}\left(T(\alpha(r+|c|), f)+\frac{|f(0)|}{2}\right)\\
&\quad + J(r+c,f)-J(r-c,f)+J(r+c,-f)-J(r-c,-f).
\end{aligned}
\]

Case 2. Suppose $c<0,$ then by the Poisson-Jensen formula in Theorem \ref{T3.1}, we have
\[
\begin{aligned}
f(&x+c)-f(x)=\frac{c}{2\rho}(f(\rho)-f(-\rho))\\
&+\frac{1}{2\rho} \sum_{|c_j|<\rho}\omega_{f}\left(c_{j}\right)\bigl((|c_{j} - x - c| - |c_{j} - x|)\rho + c c_{j}\bigr) \\
& - \frac{c}{2\rho} \Bigg( \sum_{ -\rho\leq k_j^+\leq \rho}\kern-6pt\Omega_{f}\left(k_{j}^+\right)+\sum_{-\rho\leq k_j^-\leq \rho}\kern-6pt\Omega_{f}\left(k_{j}^-\right)\Bigg)\\
&+\frac12\Bigg( \sum_{x+c< k_j^+\leq x}\kern-6pt\Omega_{f}(k_{j}^+) +\sum_{x+c< k_j^-\leq x}\kern-6pt\Omega_{f}(k_{j}^-)+\sum_{ x+c\leq k_j^+ < x}\kern-6pt\Omega_{f}(k_{j}^+)+\sum_{ x+c\leq k_j^-< x}\kern-6pt\Omega_{f}(k_{j}^-) \Bigg).
\end{aligned}
\]
For this case, we only need to estimate that
\[
\begin{aligned}
\sum_{x+c< k_j^+\leq x}\kern-6pt\Omega_{f}(k_{j}^+)+\sum_{ x+c\leq k_j^+ < x}\kern-6pt\Omega_{f}(k_{j}^+)=&\sum_{\alpha_{\mu}^+\in[x+c,x]}\kern-6pt\Omega_{f}\left(\alpha_{\mu}^+\right)+\sum_{\alpha_{\mu}^+\in(x+c,x)}\kern-6pt\Omega_{f}\left(\alpha_{\mu}^+\right)\\
&+\sum_{\beta_{\nu}^+\in[x+c,x]}\kern-6pt\Omega_{f}\left(\beta_{\nu}^+\right)+\sum_{\beta_{\nu}^+\in(x+c,x)}\kern-6pt\Omega_{f}\left(\beta_{\nu}^+\right)
\end{aligned}
\] and
\[
\begin{aligned}
\sum_{x+c< k_j^-\leq x}\kern-6pt\Omega_{f}(k_{j}^-)+\sum_{ x+c\leq k_j^- < x}\kern-6pt\Omega_{f}(k_{j}^-)=&\sum_{\xi_{s}^-\in[x+c,x]}\kern-6pt\Omega_{f}\left(\xi_{s}^-\right)+\sum_{\xi_{s}^-\in(x+c,x)}\kern-6pt\Omega_{f}\left(\xi_{s}^-\right)\\
&+\sum_{\zeta_{t}^-\in[x+c,x]}\kern-6pt\Omega_{f}\left(\zeta_{t}^-\right)+\sum_{\zeta_{t}^-\in(x+c,x)}\kern-6pt\Omega_{f}\left(\zeta_{t}^-\right).
\end{aligned}
\]
Hence, similar to the discussion for the Case 1, we have
\[
\begin{aligned}
m\biggl(r,\frac{1}{2}\biggr(\sum_{\alpha_{\mu}^+\in[x+c,x]}\kern-6pt\Omega_{f}\left(\alpha_{\mu}^+\right)&+\sum_{\alpha_{\mu}^+\in(x+c,x)}\kern-6pt\Omega_{f}\left(\alpha_{\mu}^+\right)\biggr)\biggr) \\
& \quad \le \frac12\bigl(\sum_{\alpha_{\mu}^+\in[r+c,r-c]}\kern-16pt|\Omega_f(\alpha_{\mu}^+)|
+\sum_{\alpha_{\mu}^+\in[-r+c,-r-c]}\kern-18pt|\Omega_f(\alpha_{\mu}^+)|\bigl)\\
& \quad \leq J^+(r-c,f)-J^+(r+c,f),
\end{aligned}
\]
\[
\begin{aligned}
m\biggl(r,\frac{1}{2}\biggr(\sum_{\beta_{\nu}^+\in[x+c,x]}\kern-6pt\Omega_{f}\left(\beta_{\nu}^+\right)+\sum_{\beta_{\nu}^+\in(x+c,x)}\kern-6pt\Omega_{f}\left(\beta_{\nu}^+\right)\biggr)\biggr)&\leq J^+(r-c,-f)-J^+(r+c,-f),\\
m\biggl(r,\frac{1}{2}\biggr(\sum_{\xi_{s}^-\in[x+c,x]}\kern -8pt\Omega_{f}\left(\xi_{s}^-\right)+\sum_{\xi_{s}^-\in(x+c,x)}\kern -8pt\Omega_{f}\left(\xi_{s}^-\right)\biggr)\biggr)&\leq J^-(r-c,f)-J^-(r+c,f)
\end{aligned}
\]
and
\[
\begin{aligned}
m\biggl(r,\frac{1}{2}\biggr(\sum_{\zeta_{t}^-\in[x+c,x]}\kern -8pt\Omega_{f}\left(\zeta_{t}^-\right)+\sum_{\zeta_{t}^-\in(x+c,x)}\kern -8pt\Omega_{f}\left(\zeta_{t}^-\right)\biggr)\biggr)\leq J^-(r-c,-f)-J^-(r+c,-f).
\end{aligned}
\] The remaining estimates are consistent with Case 1; so we have
\[
\begin{aligned}
m(r,f(x+c)\oslash f(x))&\leq\frac{16|c|}{r+|c|}\frac{1}{\alpha-1}\left(T(\alpha(r+|c|), f)+\frac{|f(0)|}{2}\right)\\
&\quad + J(r-c,f)-J(r+c,f)+J(r-c,-f)-J(r+c,-f).
\end{aligned}
\]

Combing Case 1 and Case 2, whenever $c>0$ or $c<0,$ we always have
\[
\begin{aligned}
m(r, f(x+c)\oslash f(x))&\leq\frac{16|c|}{r+|c|}\frac{1}{\alpha-1}\left(T(\alpha(r+|c|), f)+\frac{|f(0)|}{2}\right)\\
&\quad + J(r+|c|,f)-J(r-|c|,f)+J(r+|c|,-f)-J(r-|c|,-f).
\end{aligned}
\] Therefore, this proof is complete.
\end{proof}

Finally, we state the logarithmic derivative lemma and provide its proof.

\begin{lemma}[Logarithmic derivative lemma]\label{LDL} Let $c\in\mathbb{R}\setminus\{0\}$. If $f$ is a slope-admissible generalized tropical meromorphic function on $\mathbb{R}$ with
with hyper-order strictly less than one, then
\[m(r,f(x + c)\oslash f(x)) = o(T(r,f))\]
holds outside an exceptional set of finite logarithmic measure as $r$ tends to infinity.
\end{lemma}
\begin{proof} To apply Theorem \ref{T2.7}, define $\phi(r):=r$, $\xi(x):=(\log x)^{1+\varepsilon}$ and
\[
\alpha=1+\frac{\phi(r+|c|)}{(r+|c|)\xi(T(r+|c|,f))}
=1+\frac{1}{(\log T(r+|c|,f))^{1+\varepsilon}}.
\]
By Theorem \ref{T4.33},
\[
\begin{aligned}
m\left(r,\frac{f(x+c)}{f(x)}\oslash\right)&\leq\frac{16|c|}{r+|c|}\frac{1}{\alpha-1}\left(T(\alpha(r+|c|), f)+\frac{|f(0)|}{2}\right)\\
&\quad + J(r+|c|,f)-J(r-|c|,f)+J(r+|c|,-f)-J(r-|c|,-f)\\
&=\frac{16|c|}{r+|c|}\bigl(\log T(r+|c|,f)\bigr)^{1+\varepsilon}\left(T(\alpha(r+|c|),f)+\frac{|f(0)|}{2}\right) \\
&\quad + J(r+|c|,f)-J(r-|c|,f)+J(r+|c|,-f)-J(r-|c|,-f)
\end{aligned}
\]
hold outside of an exceptional set of finite logarithmic measure.
To proceed, fix $\varepsilon>0$ so that $(\rho_2(f)+\varepsilon)(1+\varepsilon)\le 1-\gamma$ for some $0<\gamma<1$.
Then, since $(\log T(r+|c|,f))^{1+\varepsilon}\le (r+|c|)^{1-\gamma}$,
\[
\frac{\bigl(\log T(r+|c|,f)\bigr)^{1+\varepsilon}}{r+|c|}\le (r+|c|)^{-\gamma}\to 0
\]
as $r$ approaches to infinity. For the term
\[
T(\alpha(r+|c|),f)=T\left(r+|c|+\frac{r+|c|}{(\log T(r+|c|,f))^{1+\varepsilon}},f\right)
\]
above, we may apply Theorem \ref{T2.7} to conclude that
\[
T(\alpha(r+|c|),f)\le C T(r+|c|,f)
\]
outside an exceptional set of finite logarithmic measure, and further by Theorem \ref{T4.4}, we have
\[
T(\alpha(r+|c|),f)\le C T(r,f)
\] Since $J(r,f)$ is nondecreasing and \(\limsup_{r\to \infty}\frac{\log\log J(r,f)}{\log r}<1
\) ( because $J(r,f)\leq T(r,f)$), it follows from Theorem \ref{T4.4}, with $s=2|c|$, that
\[
J(r+|c|,f)-J(r-|c|,f)=o(J(r-|c|,f))
\] as $r$ tends to infinity outside an exceptional set $E_1$ of finite logarithmic measure. Similarly,
\[
J(r+|c|,-f)-J(r-|c|,-f)=o(J(r-|c|,-f))\] as $r$ tends to infinity outside an exceptional set $E_2$ of finite logarithmic measure.
Since
\[
J(r-|c|,f)\leq T(r-|c|,f)\leq T(r,f),
\] and \[
J(r-|c|,-f)\leq T(r-|c|,-f)\leq T(r,-f)=T(r,f)-f(0),
\] we obtain \[
J(r+|c|,f)-J(r-|c|,f)+J(r+|c|,-f)-J(r-|c|,-f)=o(T(r,f))
\]
as $r$ tends to infinity outside an exceptional set $E=E_1\cup E_2$ of finite logarithmic measure.
Therefore, outside an exceptional set of finite logarithmic measure,
\[
m\left(r,\frac{f(x+c)}{f(x)}\oslash\right)
\le 32C(r+|c|)^{-\gamma}T(r,f)+o(T(r,f))=o(T(r,f)).
\]
\end{proof}

\section{Applications to ultra-discrete equations} Our main application of tropical Nevanlinna theory is to measure the complexity of solutions to ultra-discrete equations. Many such equations have been considered in the literature \cite{3, Hayman64}. Examples include
\begin{align}
y_{n+1}+y_{n-1} &= \max\{y_n+n,0\}-y_n, \label{eq:10}\\
y_{n+1}+y_{n-1} &= a+\max\{y_n,n\}-\max\{y_n+n,0\}-y_n, \label{eq:11}\\
y_{n+1}+y_{n-1} &= \max\{n+a,y_n\}+\max\{n-a,y_n\} \notag\\
&\quad -\max\{y_n+n+b,0\}-\max\{y_n+n-b,0\}, \label{eq:12}
\end{align}
where $a$ and $b$ are constants. However, the above equations only consider mappings from $\mathbb Z$ to itself. In 2009, Halburd and Southall \cite{Halburd09} extended such functions to continuous piecewise linear functions in the real variable $x$. In particular, instead of equation \eqref{eq:10}, they consider the extended equation
\begin{equation}
y(x+1)+y(x-1)=\max\{y(x)+x,0\}-y(x),\quad x\in\mathbb R. \label{eq:13}
\end{equation} In this paper, we further generalize them to possibly discontinuous piecewise-linear functions.

First, we consider the first-order ultra-discrete equations.

\begin{lemma}\label{L5.1} The equation
\begin{align}\label{eq12}
y(x+1)=y(x)^{\otimes c}=cy(x), \quad\quad c\in\mathbb{R}\setminus\{0\}
\end{align} admits a non-constant generalized tropical meromorphic solution $f$ on \(\mathbb{R}\). Moreover, this solution satisfies $\rho_2(f)<1$ if and only if $c=\pm 1$.
\end{lemma}
\begin{proof} We first prove the sufficiency. If $c=1$, then $y(x)$ is an arbitrary one-periodic function. Clearly, it suffices to restrict this one-periodic function to a generalized tropical meromorphic function of finite order, and the conclusion follows. If $c=-1$, then
\[
y(x+1)=-y(x),\quad y(x+2)=-y(x+1).
\] Hence, we see that $y(x+2)=y(x)$, that is, $y(x)$ is arbitrary one-antiperiodic function. Similarly, equation \eqref{eq12} also admits a generalized tropical meromorphic solution of finite order.\par

We now prove the necessity. Since $y(x)$ is non-constant, it must have either a point of discontinuity in its derivative or a jump discontinuity on $[0,1)$.\par

Suppose there exists an $x_0\in [0,1)$, such that $x_0$ is a pole of $y(x)$, then we have
\begin{align*}
N(r,f)&=\frac{1}{2}\sum_{|b_v|<r}|\omega_f(b_v)|(r - |b_{\nu}|)\\
&\geq \frac{1}{2}\bigl(c^{[r-2]}|\omega_f(x_0)|(2-x_0)+c^{[-r+1]}|\omega_f(x_0)|(1+x_0)\bigl).
\end{align*}
Assume that $c\neq \pm 1$, then we can see that
\[
\rho_2(f)=\limsup_{r\to \infty}\frac{\log\log T(r,f)}{\log r}\geq \limsup_{r\to \infty}\frac{\log\log N(r,f)}{\log r}=1.
\] This leads to a contradiction with $\rho_2(f)<1$.\\
Suppose there exists an $x_0\in [0,1)$, such that $x_0$ is a root of $y(x)$, by $$T(r, f)=T(r, 1_{\mathbb{T}}\oslash f)+f(0),$$  we similarly obtain
\[
\rho_2(f)=\limsup_{r\to \infty}\frac{\log\log T(r,f)}{\log r}=\limsup_{r\to \infty}\frac{\log\log \bigl(T(r,1_{\mathbb{T}}\oslash f)+f(0)\bigl)}{\log r}\geq 1.
\] again contradicting $\rho_2(f)<1$.\par

Suppose there exists an $x_0\in [0,1)$, such that $x_0$ is a right-negative jump of $y(x)$, then we have
\begin{align*}
J^+(r,f)&=\frac{1}{2} \sum_{|\beta_\nu^+| \leq r} h_f(\beta_\nu^+)\\
&\geq \frac{1}{2}\bigl(c^{[r-1]}+c^{[-r]}|\bigl)|\Omega_f(x_0)|.
\end{align*}
Assume that $c\neq \pm 1$, then we can see that
\[
\rho_2(f)=\limsup_{r\to \infty}\frac{\log\log T(r,f)}{\log r}\geq \limsup_{r\to \infty}\frac{\log\log J^+(r,f)}{\log r}=1.
\] This also leads to a contradiction with $\rho_2(f)<1$.\\
If $x_0\in [0,1)$ is a right-positive jump, a left-negative jump, or a left-positive jump point, the same argument as above yields $\rho_2(f)\geq 1$.
\end{proof}

Next, we consider the following second-order ultra-discrete equations.
\begin{lemma}\label{L5.2}
Let $y\not\equiv 0$ be a generalized tropical meromorphic solution of
\begin{align}\label{eq13}
y(x+1)\otimes y(x-1)=y(x)^{\otimes c},
\end{align}  for some  $c\in\mathbb{R}$. If $y$ is of finite order then $|c|\le 2$.
\end{lemma}
\begin{proof} Assume that \(|c|> 2\). Note that the equation \eqref{eq13} can be rewritten as
\begin{align}\label{eq14}
y(x+1)-ay(x)=b(y(x)-ay(x-1)),
\end{align}
where $a$ and $b$ are the roots of the quadratic equation $\lambda^2-c\lambda+1=0$, that is, $a+b=c$ and $ab=1$. Without loss of generality, set
\[
a=\frac{c+\sqrt{c^2-4}}{2},~\quad b=\frac{c-\sqrt{c^2-4}}{2}.
\]

Since $y\not\equiv 0$, there exists $x_0\in[0,2)$ such that $y(x_0)\neq 0$. Therefore, for each $\nu\in\mathbb Z$,
\begin{equation*}
y(x_0+\nu)=\alpha a^{\nu}+\beta b^{\nu},
\end{equation*}
where
\[
\alpha=\frac{y(x_0+1)-by(x_0)}{a-b}
\quad\text{and}\quad
\beta=\frac{y(x_0+1)-ay(x_0)}{b-a}
\]
are not both zero.

Now if $c>2$, then $a>1$ and $b^{-1}>1$, while if $c<-2$ then $-b>1$ and $-a^{-1}>1$. By an argument analogous to that in Lemma \ref{L5.1}, we obtain $\rho_2(y)\geq1$. Hence $y(x)$ has infinite order, which leads to a contradiction.
\end{proof}

In light of Lemma \ref{L5.2}, we now consider the following ultra-discrete equation.
\begin{equation}
y(x+1)\otimes y(x-1)=R\big(x,y(x)\big),
\label{eq:15}
\end{equation} where $R$ is tropical rational in $x$ and $y$. We remark that all such equations admit infinitely many generalized tropical meromorphic solutions.

\begin{theorem}\label{T5.3}
Let $P(y)=\max\{a_0,a_1+y,\dots,a_p+py\}$ and $Q(y)=\max\{b_0,b_1+y,\dots,b_q+qy\}$ be two tropical polynomials with no common roots and neither $a_p$ nor $b_q$ is $-\infty$. If $|p-q|>2$, then the equation
\begin{equation}\label{eq15}
y(x+1)+y(x-1)=P(y(x))-Q(y(x))
\end{equation}
admits infinitely many infinite-order generalized tropical meromorphic solutions.
\end{theorem}

\begin{proof}
If $y$ is sufficiently large for all $x$ larger than some number $\xi$, then equation \eqref{eq15} reduces to
\[
y(x+1)+y(x-1)=(p-q)y(x)+a_p-b_q,
\]
for all $x>\xi$. By the proof of \cite[Theorem 4.3]{Halburd09}, we know that if $p-q>2$ then, given $\alpha>0$, there is a family of solutions such that any member evaluated at an integer $\nu>\xi$ has the form
\[
y(\nu)=\frac{b_q-a_p}{p-q-2}+\alpha\left(\frac{(p-q)+\sqrt{(p-q)^2-4}}{2}\right)^{\nu}.
\] By an argument analogous to that in Lemma \ref{L5.1}, we obtain $\rho_2(y)\geq 1$. Hence $y(x)$ has infinite order. If $p-q<-2$ then the same argument works with a minus sign in front of the square root.
\end{proof}

Lastly, by means of the logarithmic derivative lemma  \ref{LDL} and \cite[Theorem 4.5]{Halburd09}, we will present a Clunie's lemma for slope-admissible piecewise-linear functions. Let $\lambda=(\lambda_0,\lambda_1,\dots,\lambda_m)$ be a multi-index of non-negative real numbers, and consider
\begin{align*}
P_{\lambda}(x,f)&:=\bigotimes^m_{j=0}f(x+c_j)^{\otimes \lambda_j}\\
&=\lambda_0 f(x)+\lambda_1 f(x+c_1)+\cdots+\lambda_m f(x+c_m)
\end{align*}
with given real shifts $(0,c_1,\dots,c_m)\in\mathbb{R}^{m+1}$. Then, an expression of the form
\[
P(x,f)=\bigoplus_{\lambda\in\Lambda}a_\lambda(x)\otimes P_{\lambda}(x,f)= \bigoplus_{\lambda\in\Lambda} a_\lambda(x)\otimes \bigotimes^m_{j=0}f(x+c_j)^{\otimes \lambda_j},
\] where $\Lambda$ is a finite set of indices, is called a tropical difference polynomial of total degree
\[
\deg(P):=\max_{\lambda\in\Lambda[P]}||\lambda||, \quad ||\lambda||:=\lambda_0+\dots+\lambda_m
\]
in $f$ and its shifts. We will say that the coefficients are small if $T(r,a_\lambda)=o\bigl(T(r,f)\bigr)$ outside a set of finite logarithmic measure.

\begin{theorem}\label{T5.4}
Let $P(x,f)$ and $Q(x,f)$ be tropical difference polynomials in $f$ and its shifts with small coefficients. If $f$ is a slope-admissible generalized tropical meromorphic function satisfying hyper-order strictly less than one and $f^{\otimes n}(x)P(x,f)=Q(x,f)$, where the degree of $Q$ in $f$ and its shifts is less than or equal to $n$, then
\[
m(r,P(x,f))=o\left(T(r,f)\right)
\]
holds outside an exceptional set of finite logarithmic measure as $r$ tends to infinity.
\end{theorem}

\begin{proof}
Given a fixed $r>0$, let
$$S_+:=\{s:\,f(s)\ge0\text{ and }|s|=r\} \quad \text{and}\quad S_-:=\{s:\,f(s)<0\text{ and }|s|=r\}.$$ Then
\[
m(r,P(x,f))=\frac12\bigl(\sum_{s\in S_+}P(s,f)^+ + \sum_{s\in S_-}P(s,f)^+\bigl).
\]
We put
$$P(x,f)=\bigoplus_{\lambda\in\Lambda_P} a_\lambda(x)\otimes f^{\otimes\lambda}(x)\quad \text{and}\quad Q(x,f)=\bigoplus_{\lambda\in\Lambda_Q} b_\lambda(x)\otimes f^{\otimes\lambda}(x).$$

Suppose that $s\in S_-$, then by $f(s)<0$, we have
\begin{align*}
P(s,f)&=\bigoplus_{\lambda\in\Lambda_P} a_\lambda(s)\otimes f^{\otimes\lambda_0}(s)\otimes f^{\otimes\lambda_1}(s+c_1)\otimes\cdots\otimes f^{\otimes\lambda_m}(s+c_m)\\
&\le \max_{\lambda\in\Lambda_P}
\Bigl\{a_\lambda(s)+\lambda_1\bigl[f(s+c_1)-f(s)\bigr]+\cdots+\lambda_m\bigl[f(s+c_m)-f(s)\bigr]\Bigr\}\\
&\le \max_{\lambda\in\Lambda_P}|a_\lambda(s)|+\max_{\lambda\in\Lambda_P}\sum^m_{j=0}\lambda_j|f(s+c_j)-f(s)|\\
&\le 2\max_{\lambda\in\Lambda_P}m(r,|a_\lambda|)+2\max_{\lambda\in\Lambda_P}\sum^m_{j=0}\lambda_j m(r,|f(x + c_j)\oslash f(x)|) .
\end{align*} Since \[
|f(x + c_j)\oslash f(x)|=\bigl(f(x + c_j)\oslash f(x)\bigl) \oplus \bigl(-f(x + c_j)\oslash -f(x)\bigl).
\] it follows from Lemma \ref{L2.3}, Theorem \ref{T3.2} and Lemma \ref{LDL} that \begin{align*}
m(r,|f(x + c_j)\oslash f(x)|)&=m\bigl(r,\bigl(f(x + c_j)\oslash f(x)\bigl) \oplus \bigl(-f(x + c_j)\oslash -f(x)\bigl)\bigl)\\
&\leq m\bigl(r,f(x + c_j)\oslash f(x)\bigl)+ m\bigl(r,\bigl(-f(x + c_j)\oslash -f(x)\bigl)\bigl)\\
&=o\left(T(r,f)\right)+o\left(T(r,-f)\right)\\
&=o\left(T(r,f)\right).
\end{align*}
Hence, we have
\[
\sum_{s\in S_-} P(s,f)^+ = o\left(T(r,f)\right),
\]
outside an exceptional set of finite logarithmic measure.

Suppose that $s\in S_+$, we note that $P(s,f)=Q(s,f)-nf$ and the degree of $Q$ is at most $n$. Hence, we similarly obtain
\begin{align*}
P(s,f)&\le \max_{\lambda\in\Lambda_Q}
\Bigl\{b_\lambda(s)+\lambda_1\bigl[f(s+c_1)-f(s)\bigr]+\cdots+\lambda_m\bigl[f(s+c_m)-f(s)\bigr]\Bigr\}\\
&\le 2\max_{\lambda\in\Lambda_Q}m(r,|b_\lambda|)+2\max_{\lambda\in\Lambda_Q}\sum^m_{j=0}\lambda_j m(r,|f(x + c_j)\oslash f(x)|).
\end{align*}
So again using Lemma \ref{LDL}, we have
\[
\sum_{s\in S_+} P(s,f)^+ = o\left(T(r,f)\right),
\]
outside an exceptional set of finite logarithmic measure.
\end{proof}


\begin{thebibliography}{99}
\bibitem{1} T. M. Apostol, Mathematical Analysis, second ed., Addison-Wesley Publishing Co., Reading, Mass.-London-Don Mills, Ont., 1974. MR 0344384
\bibitem{Cao-Zheng25} T. B. Cao, J. H. Zheng, Nevanlinna theory for tropical hypersurfaces, Annali della Scuola Normale Superiore di Pisa, Classe di Scienze 16(5), 2025: 2105-2144.
\bibitem{Cao-peng} T. B. Cao and J. H. Peng, Tropical Nevanlinna theory in several variables, Journal of the London Mathematical Society 113(2), 2026: e70449, pp1-49, DOI: 10.1112/jlms.70449
\bibitem{3} B. Grammaticos, Y. Ohta, A. Ramani, D. Takahashi, and K. M. Tamizhmani. Cellular automata and ultra-discrete Painlev'e equations. Phys. Lett. A, 226:53–58, 1997.
\bibitem{Halonen23} J. Halonen, R. J. Korhonen and G. Filipuk, Tropical second main theorem and the Nevanlinna inverse problem, arXiv:2305.13939, (2023).
\bibitem{Halonen24} J. Halonen, A tropical analogue of the lemma on the logarithmic derivative, Bulletin of the Australian Mathematical Society 110.2 (2024): 303-312.
\bibitem{Halburd09} R. G. Halburd, N. J. Southall, Tropical Nevanlinna theory and ultradiscrete equations, International Mathematics Research Notices 2009.5 (2009): 887-911.
\bibitem{Hayman64} W. K. Hayman. Meromorphic functions. Oxford Mathematical Monographs. Clarendon Press, Oxford, 1964.
\bibitem{Risto15} R. J. Korhonen, I. Laine, K. Tohge, Tropical value distribution theory and ultra-discrete equations, World Scientific Co. Pte. Ltd., Singapore, 2015.
\bibitem{Speyer09} D. Speyer, B. Sturmfels, Tropical mathematics, Mathematics Magazine 2009, 82 (3), 163-173.
\bibitem{Risto26}R. J. Korhonen and C. L. Tan, $n$-th tropical Nevanlinna theory, arXiv:2602.03500, (2026)
\bibitem{Risto16} R. J. Korhonen and K. Tohge, Second main theorem in the tropical projective space, Advances in Mathematics 298(2016): 693-725.
\bibitem{Laine11} I. Laine and K. Tohge, Tropical Nevanlinna theory and second main theorem, Proceedings of the London Mathematical Society 102(5), 2011: 883-922.
\bibitem{mik} G. Mikhalkin, Enumerative tropical algebraic geometry in $\mathbb{R}^2,$ J. Amer. Math. Soc. 18(2005), No.2, 313-377.
\bibitem{Maclagan-15} D. Maclagan, B. Sturmfels, Introduction to tropical geometry, vol. 161, American Mathematical Soc, RI, 2015.
\bibitem{W-C26} Y. T. Wang, T. B. Cao, Tropical Cartan's second main theorem for hyperplanes in general position. arXiv preprint arXiv:2605.23687,  (2026).
\bibitem{Mik2006} G. Mikhalkin, Tropical geometry and its applications, International Congress of Mathematicians, vol. II,
European Mathematical Society, Zürich, 2006, pp. 827-852.
\bibitem{Mik2020} G. Mikhalkin and J. Rau, Tropical Geometry, IRMA Lectures in Mathematics and Theoretical Physics, vol. 35,
European Mathematical Society, Zürich, 2020.
\bibitem{Orm2010} C. M. Ormerod,  Hypergeometric solutions to an ultradiscrete Painlevé equation, Journal of Nonlinear Mathematical Physics 17(1), 2010: 87-102
\bibitem{Laine16} I. Laine, K. Liu, K. Tohge, Tropical variants of some complex analysis results, Annales Fennici Mathematici 41(2), 2016: 923-946.
\bibitem{Zhang18} L. Zhang,  G. Naitzat, and L. H. Lim.  Tropical geometry of deep neural networks.  International Conference on Machine Learning. PMLR, 2018.

\end{thebibliography}
\end{document}